\documentclass[10pt,reqno,onecolumn]{elsart}
\usepackage{amsmath,amssymb,color}
\usepackage[T1]{fontenc}
\usepackage{natbib}
\usepackage{graphicx}
\usepackage{epstopdf}
\usepackage{subfigure}
\usepackage{rotating}

\def\dd{\mathrm d}
\def\bh{\mathbf h}
\def\bb{\mathbf b}

\def\bD{\mathbf D}

\def\bH{\mathbf H}
\def\bQ{\mathbf Q}
\def\bh{\mathbf h}

\def\Sset{\mathbb S}
\def\Hset{\mathbb H}
\def\Lset{\mathbb L}
\def\Nset{\mathbb N}
\def\Rset{\mathbb R}

\def\Dset{\mathbb D}
\def\lm{_{\ell m}}
\def\tr{\text{tr }}
\def\ie{\emph{i.e.}}
\def\eg{\emph{e.g.}}
\def\dint{\int\!\!\!\int}

\def\esp{\mathbb E}
\def\proba{\mathbb P}
\def\lmin{\ell_{\min}}
\def\lmax{\ell_{\max}}

\def\ind{{\mathbf{1}}} 
\def\sumstar{\sum\nolimits^*}

\newcommand\wigner[6]{\textstyle\left( \scriptstyle\begin{array}{ccc}
    #1 & #2 & #3\\
    #4 & #5 & #6
  \end{array}
  \right)}

\newcommand\varwigner[6]{{\textstyle\left[
      \scriptstyle\begin{array}{ccc}
    #1 & #2 & #3\\
    #4 & #5 & #6
  \end{array}
  \right]}}

\def\mcJ{\mathcal J}

 {\everymath{\displaystyle\everymath{}}\array}%
 {\endarray}

\newlength\jataille
\newcommand{\figgauche}[3]%
{\jataille=\\advance\jataille by -#1
\advance\jataille by -.5cm
\begin{minipage}[c]{#1}
   \includegraphics[width=#1]{#2}
\end{minipage}\hfill
\begin{minipage}[c]{\jataille}
   \footnotesize #3 \normalsize
\end{minipage}}

\numberwithin{thm}{section}

\newtheorem{hyp}[thm]{Assumption}

\begin{document}

\begin{frontmatter}
  \title{Practical wavelet design on the sphere}
  \author[APC,LPMA,LTCI]{Fr{\'e}d{\'e}ric Guilloux},
    \author[APC,Lille]{Gilles Fa{\"y}\corauthref{corresp}} 
    \ead{gilles.fay@univ-lille1.fr}
and
  \author[APC,LTCI]{Jean-Fran\c cois Cardoso}.  \date{\today}
  
  \address[APC]{AstroParticule et Cosmologie, CNRS and
    Universit\'e Paris 7.}  \address[LPMA]{Laboratoire Probabilit\'es
    et Mod\`eles Al\'eatoires, CNRS and Universit\'es Paris 6-7.}
  \address[LTCI]{Laboratoire du Traitement et de la Communication de
    l'Information, CNRS and T\'el\'ecom Paris.}
  \address[Lille]{Laboratoire Paul-Painlev\'e, Universit\'e Lille-1}

  \corauth[corresp]{Corresponding author.}
  
  \begin{abstract}
    We address the question of designing isotropic analysis
    functions on the sphere which are perfectly limited in the spectral domain
    and optimally localized in the spatial domain.  This work is motivated by
    the need of localized analysis tools in domains where the data is lying on
    the sphere, \eg{} the science of the Cosmic Microwave Background.  Our
    construction is derived from the localized frames introduced by
    \cite{narcowich:petrushev:ward:2006}.  The analysis frames are optimized
    for given applications and compared numerically using various criteria.
  \end{abstract}

  \maketitle

\end{frontmatter}

\newcommand\jf[1]{---\texttt{#1}---}

\section*{Introduction}

Localized analysis for spherical data has motivated many researches
during the past decade. Data defined on the sphere are studied in
domains as various as cosmology
\citep{Hinshaw+2006,Hivon+2002,McEwen+2007}, geophysics
\citep{Holschneider+2003,Wieczorek+2005}, medicine, computer
vision. When dealing with data on the whole sphere, spectral analysis
can be achieved by Spherical Harmonics Transform (SHT) -- the
equivalent of the Fourier Series on the circle. But in many practical
situations, data are defined or available on a subset of the sphere
only. For example, cosmologists try to give sharp estimates of the
cosmic microwave background (CMB) or its power spectrum but strong
foreground emissions superimpose to the CMB making the observations
unreliable for CMB studies. Moreover, fully observed clean
non stationary fields or stationary fields with additive non-stationary
noise still require spatially localized tools.  In such situations,
the SHT is not adequate, because of the poor spatial localization of
the basis functions. In the case of Euclidean spaces, in which the
Fourier Transform suffer from the same lack of localization,
multiscale and wavelets theory provide a mathematically elegant
solution of proven practical efficiency.

Adaptation to the sphere of the ``wavelet'' transform (in the broad sense of
filtering by spatially and spectrally localized functions) was introduced a
dozen years ago \citep{Schroder:Sweldens:95,Torresani95,
  Dahlke+95,Narcowich+96,Potts:Tasche:95,Freeden+97}.  Since then, Antoine \&
Vandergheynst (1999)\nocite{antoine:vandergheynst:1999} showed that any
Continuous Wavelet Transform (CWT) on the sphere can be viewed locally as a
regular CWT on the Euclidean tangent planes, thanks to the stereographic
correspondence between the sphere and the plane
\citep{antoine:vandergheynst:1999,Wiaux+2005}. One can then ``forget'' the
sphere by projecting it on tangent planes, realizing the analysis in these
planes, and then apply the inverse projection to get back eventually to the
sphere. A discretized version of this approach of CWT has been presented by
\cite{bogdanova_etal:2005}, leading to wavelet frames.
This approach has already been followed in astrophysics for the
analysis of the Cosmic Microwave Background (CMB)
\citep{Vielva+2004,McEwen+2007}. 
However these wavelets are usually defined in the spatial domain and have
infinite support in the frequency domain (which must be truncated in practice).

In the present work, we follow and extend the approach of
\cite{narcowich:petrushev:ward:2006} and their construction of ``needlets''. A
similar construction can be found in \cite{Starck+2006}.  The needlet transform
has important characteristics. Firstly it is intrinsically spherical. No
intermediate tangent plane is needed to define it. Secondly, it does not depend
on the particular spherical pixelization chosen to describe the data. Thirdly,
although the needlets still have an excellent spatial localization, they have a
finite spectral support adjustable at will . They are axisymmetric (which is
convenient when dealing with statistically isotropic random fields) and thus
the needlet coefficients are easily computed in the Spherical Harmonics
(Fourier) domain. Data filtering is defined by multiplication of the Spherical
Harmonics coefficients by well chosen window functions (which is equivalent to
convolution in spatial domain).
Needlets are well defined in theory and the statistical properties of
their coefficients have already been established for isotropic
Gaussian fields~(\cite{baldi:etal:2006a}).  However, the performance
of a needlet-based analysis depends on the particular shape of the
needlet.

This paper considers spherical filters which are generalizations of needlets in
the spirit of dual (non-tight) analysis and reconstruction frames.  We focus on
the design issue, namely the optimization of the window functions (that define
the isotropic filtering operations) for some given tasks.
We consider only band-limited needlets.  This choice is motivated by
applications in high-precision cosmology.  Indeed, the CMB power spectrum is
highly dynamic (few peaks and power-law decay) and good subsequent
cosmological parameters estimation requires high accuracy in some critically
delimited spectral ranges. Once the range is fixed, we optimize the shape of
window functions in two directions: 1) By requesting the best spatial
localization of associated needlets, in an energy-sense ($\mathbb L^2$) which
is easily solved.  This is an application of the work of \cite{Simons+2006}
which adapted to the sphere the problem solved by \cite{Slepian+78} on the real
line, giving rise to the well known prolate spheroidal wave functions (PSWF).
2) By following statistical considerations: given some region (``mask'') in
which the data is missing or thrown away and assuming that the full data is the
realization of some Gaussian isotropic random field (this is the usual
assumption made on the CMB), we minimize the mean integrated square error due
to the mask in the needlet analysis outside the mask. More criteria and
applications to cosmological science will be given in a future work.

The paper is organised as follows. In Section~\ref{sec:needlets}, we expose the
general construction of needlets. In Section~\ref{sec:crit-filt-design}, we
define and optimize the two criteria (geometrical and statistical) which
provide localized analysis filters. Their efficiency is illustrated in
Section~\ref{sec:results} with numerical simulations following the model of a
masked observation of the CMB. The proofs are postponed to
Appendix~\ref{sec:proofs}.

\section{Needlets frames}
\label{sec:needlets}

\subsection{Background and notations}

Denote $\Sset$ the unit sphere in $\Rset^3$ with generic element $\xi
= (\theta,\varphi)$ in spherical polar coordinates:
$\theta\in[0,\pi]$ is the colatitude and $\varphi\in[0,2\pi[$ the
longitude. Let $\Hset=\Lset^2(\Sset)$ be the space of complex-valued
square integrable functions on $\Sset$ under the Lebesgue measure $\dd
\xi=\sin\theta \dd \theta \dd \varphi$. Endowed with
the inner product $\langle f , g \rangle := \int_\Sset f(\xi)g^*(\xi)
\dd\xi$, $\Hset$ is a Hilbert space. Let $\| \cdot \|$ denote the
associated norm on $\Hset$.  The usual complex spherical harmonics on
$\Sset$ (which definition is recalled in Appendix \ref{sec:leg-pol})
are denoted $(Y\lm)_{\ell \geq 0, -\ell \leq m \leq \ell}$. They form
an orthonormal basis of $\Hset$.

In the following, we consider a field $X \in \Hset$. Its random spherical
harmonics coefficients or multipole moments are denoted $a\lm =\langle
X,Y\lm\rangle$. $\Hset$ can be decomposed in harmonic subspaces:
$\Hset=\bigoplus\limits^\perp_{\ell\geq0} \Hset_\ell$, where $\Hset_\ell$ is
the linear span of $Y_{\ell m} , m=-\ell,\cdots,\ell$. The number $\ell$ is
referred to as the multipole number or frequency (understood as a spatial
frequency). Let $\Pi_\ell$ be the orthogonal projection on $\Hset_\ell$. It has
an expression involving Legendre polynomials $L_\ell$ (see Appendix
\ref{sec:leg-pol})
\begin{equation}
  \label{eq:kernel}
  \Pi_\ell X(\xi)  =  \sum_{m = -\ell}^{\ell}\langle X,Y\lm\rangle Y\lm(\xi) =
  \int_\Sset X(\xi')L_\ell( \xi\cdot \xi') \dd \xi'.   
\end{equation}
where $\xi\cdot \xi' = \cos\theta\cos\theta' + \sin \theta \sin
\theta' \cos(\varphi - \varphi')$ is the usual dot product on $\Sset$.

\begin{figure}[htbp]
\centering\includegraphics{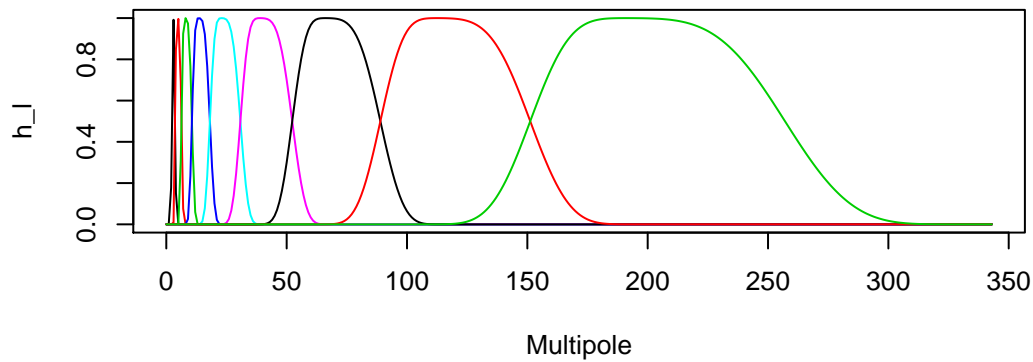}
\caption{First 10 windows satisfying conditions (\ref{reconstruction})
  and (\ref{eq:Badic}), $\mathsf h$ being a spline of order 7
  compactly supported on $[\frac1B,B]$ with
  $B=1.7$.\label{fig:filt}}
\end{figure}

A mapping on $\Sset$ which depends on the colatitude $\theta$ only is
said to be axisymmetric. The convolution of a bounded axisymmetric
function $H(\xi) = h(\cos\theta)$ with an arbitrary spherical function $X$ is
well defined through
\begin{equation}
\label{eq:defconvolution}
  H * X (\xi) = \int_\Sset h(\xi \cdot \xi') X(\xi') \dd \xi'
\end{equation}
The convolution theorem holds:
\begin{equation}
\label{eq:convolution}
  H * X = \sum\lm  h_\ell a\lm Y\lm.
\end{equation}
where $a\lm=\langle X,Y\lm\rangle$ are the  multipole moments of
$X$ and $h_\ell$ are the Legendre series coefficients  of $h$, \ie $h = \sum_{\ell \in \Nset} h_\ell L_\ell$.
Then, an isotropic wavelet analysis can be implemented either in the spatial
(\ie~direct) domain using~(\ref{eq:defconvolution}) or in the harmonic domain
using~(\ref{eq:convolution}). We choose the latter,
which accounts to multiply the harmonic coefficients of the field of
interest $X$ by a spectral window $(h_\ell)$.  
For a countable index set $\mcJ$, let
$(h^{(j)})_{j\in\mcJ}$ be a family of window functions in harmonic
domain~: $h^{(j)}\in\ell^\infty(\Nset)$. The corresponding harmonic
smoothing operators on $\Hset$ are
\begin{equation}
\label{eq:defsmoothing}
\Psi^{(j)} 
=
 \sum_{\ell \in \Nset} h^{(j)}_\ell \Pi_\ell.
\end{equation}
We call \emph{exact reconstruction condition} the one ensuring that
$\sum\limits_{j \in \mcJ} \Psi^{(j)}
=
\mathbf{Id}.$
 It also writes
\begin{equation}
\sum_{j \in \mcJ} h^{(j)}
\equiv
1\label{reconstruction}
\end{equation}
In the following, $j$ is referred to as the \emph{scale}, in analogy with the
multiresolution analysis terminology.  Important examples of windows families
having the property~(\ref{reconstruction}) may be obtained thanks to the $B$-adic mechanism: let $B >
1$, $\mcJ=\{-1\} \cup \Nset$, $h^{(-1)}_\ell = \delta_0(\ell)$ and the spectral
windows be all related to a continuous function $\mathsf h$ by
\begin{equation}
  \label{eq:Badic}
  \forall j \in \Nset, \; 
  h^{(j)}_\ell  =   \mathsf h\left(\frac\ell{B^j}\right).
\end{equation}
If $\mathsf h$ is compactly supported on $[\frac1B,B]$, then each
window $h^{(j)}$ may overlap with adjacent windows $h^{(j-1)}$ and
$h^{(j+1)}$ only. The exact reconstruction condition in this case is satisfied
as soon as
\begin{equation}
  \forall x\in[1,B],\ \mathsf h(x)+\mathsf h(B^{-1}x)=1
\end{equation}
This example is illustrated in Figures~\ref{fig:filt}
and~\ref{fig:sph-liss}.

\newlength{\widthsky}
\setlength{\widthsky}{3.5cm}
\begin{figure}[htbp]
\centering    
\subfigure[Original map]{
  \includegraphics[width=5cm]{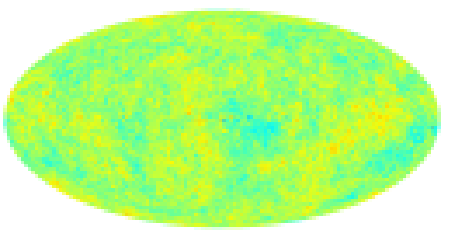}}
\subfigure[Smoothed maps, scales $j=2,...,5$]{
    \includegraphics[width=\widthsky]{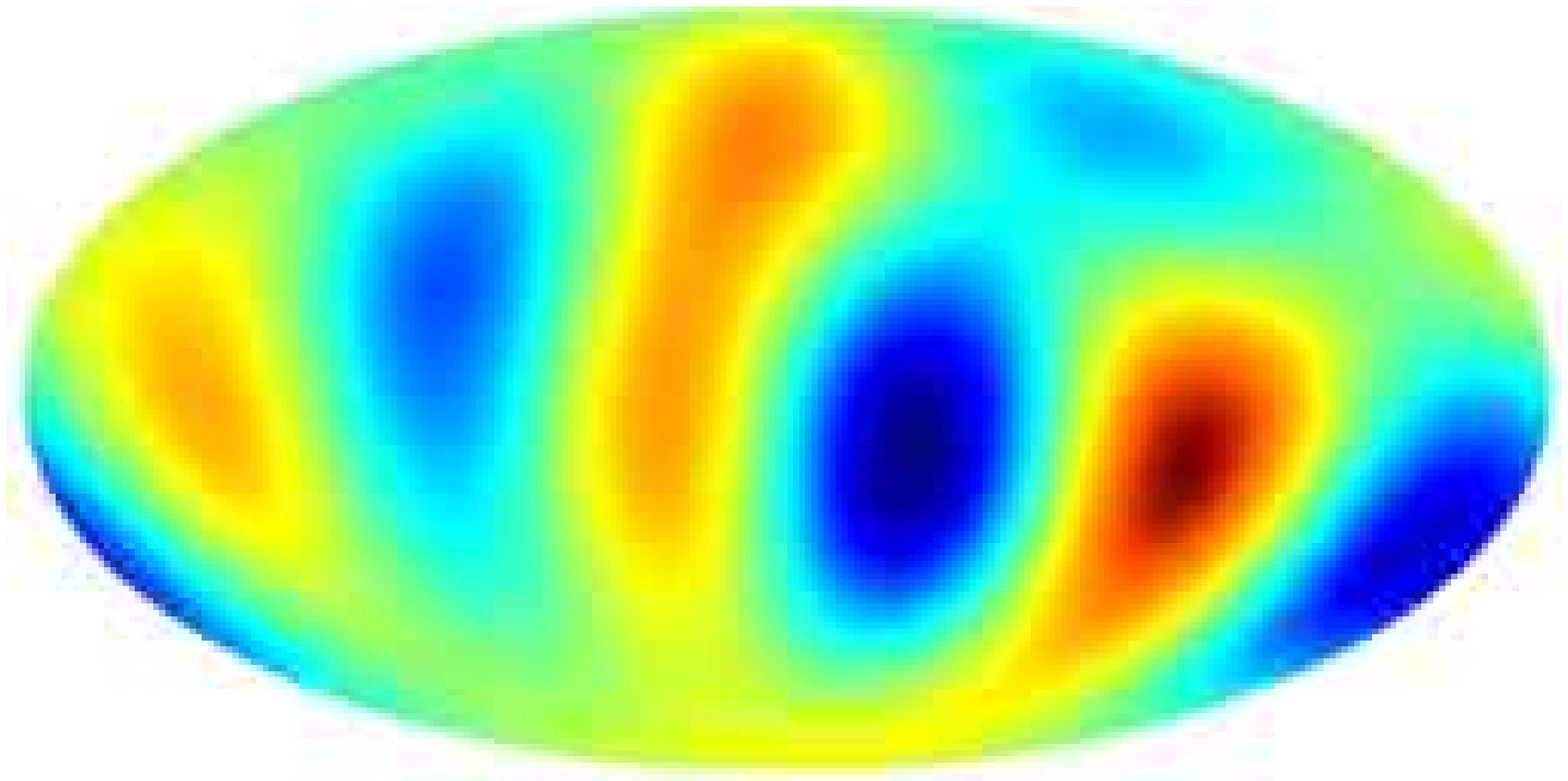}
    \includegraphics[width=\widthsky]{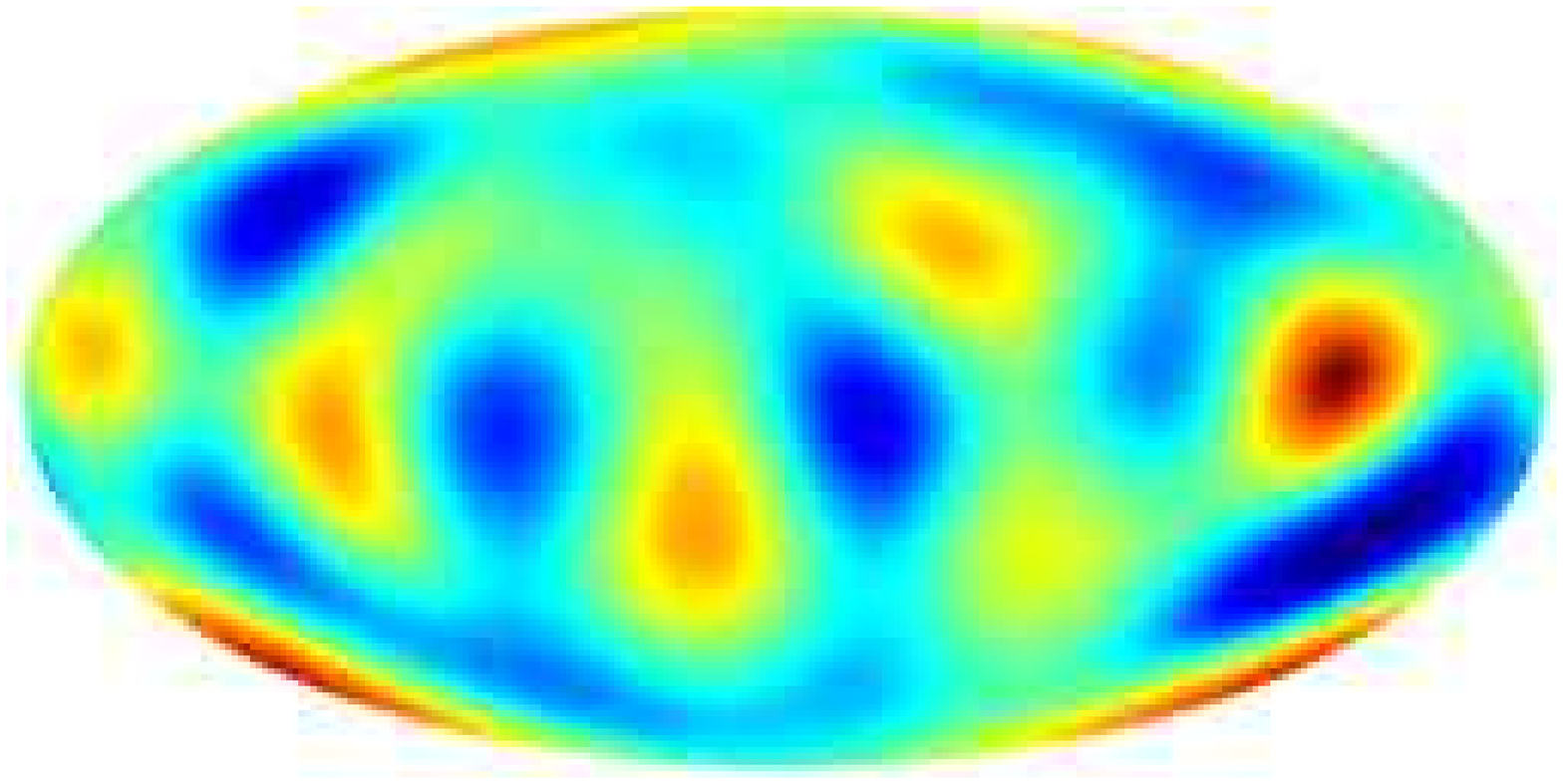}
    \includegraphics[width=\widthsky]{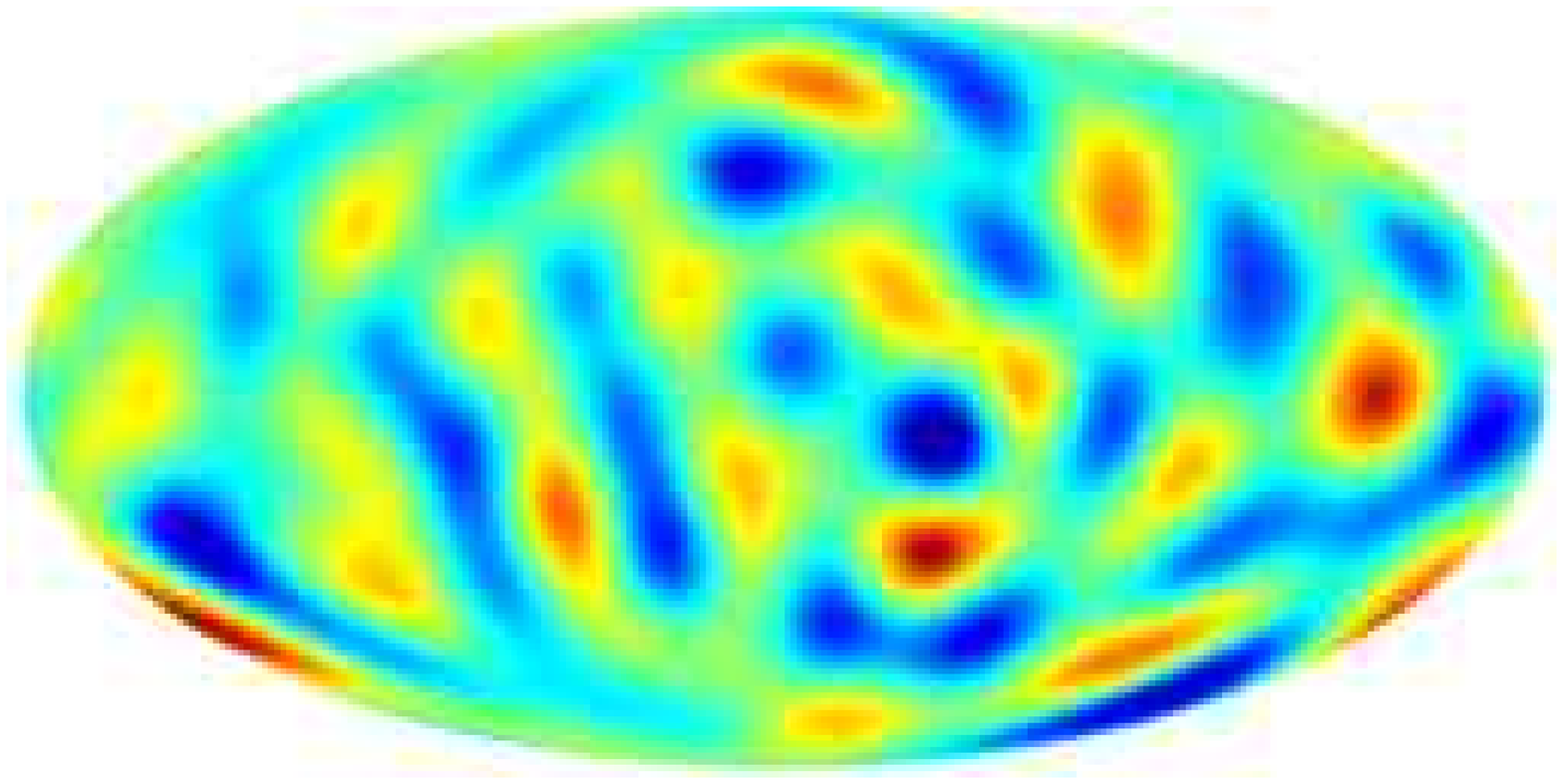}
    \includegraphics[width=\widthsky]{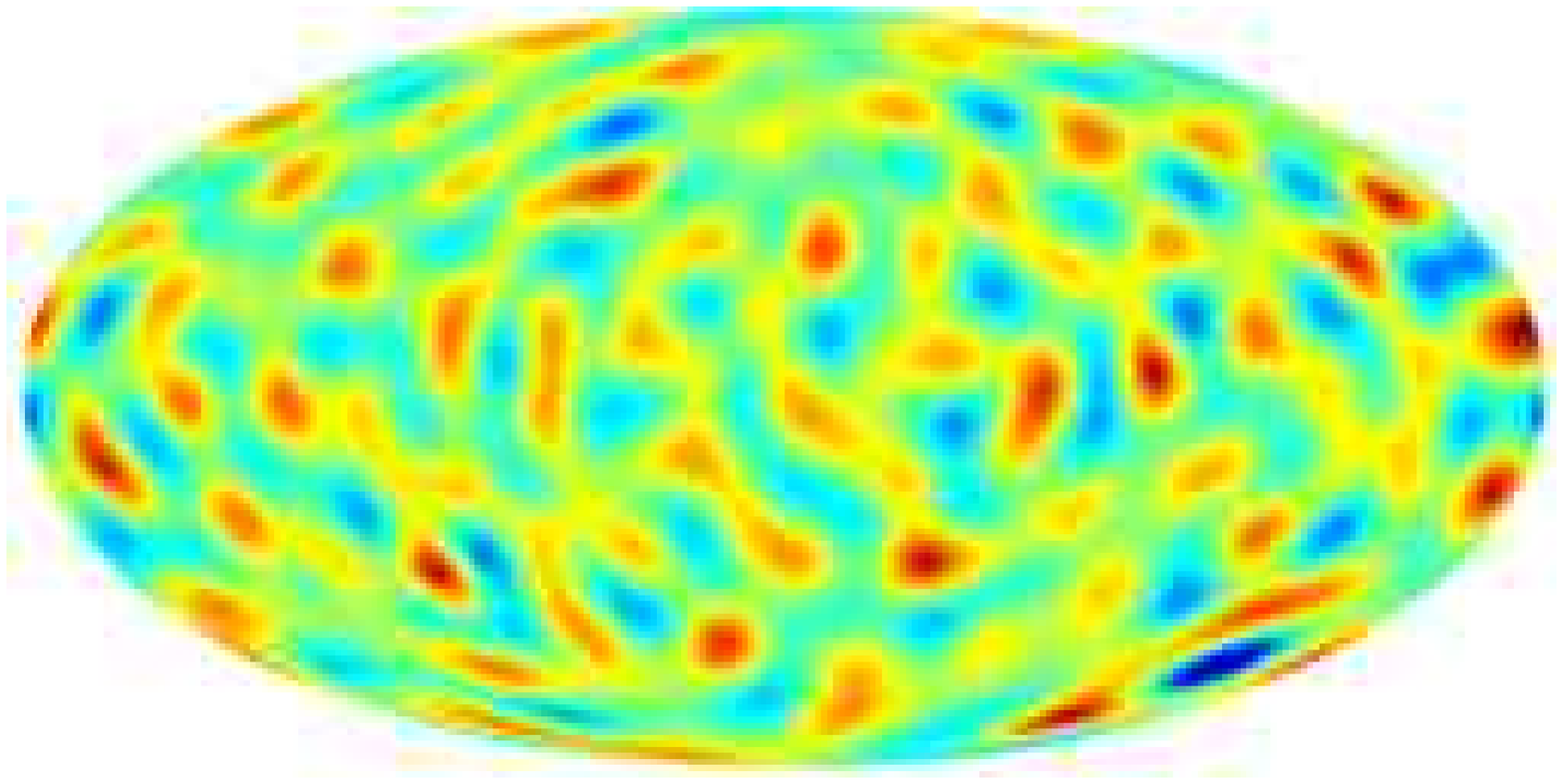}}
\subfigure[Smoothed maps, scales $j=6,...,9$]{
    \includegraphics[width=\widthsky]{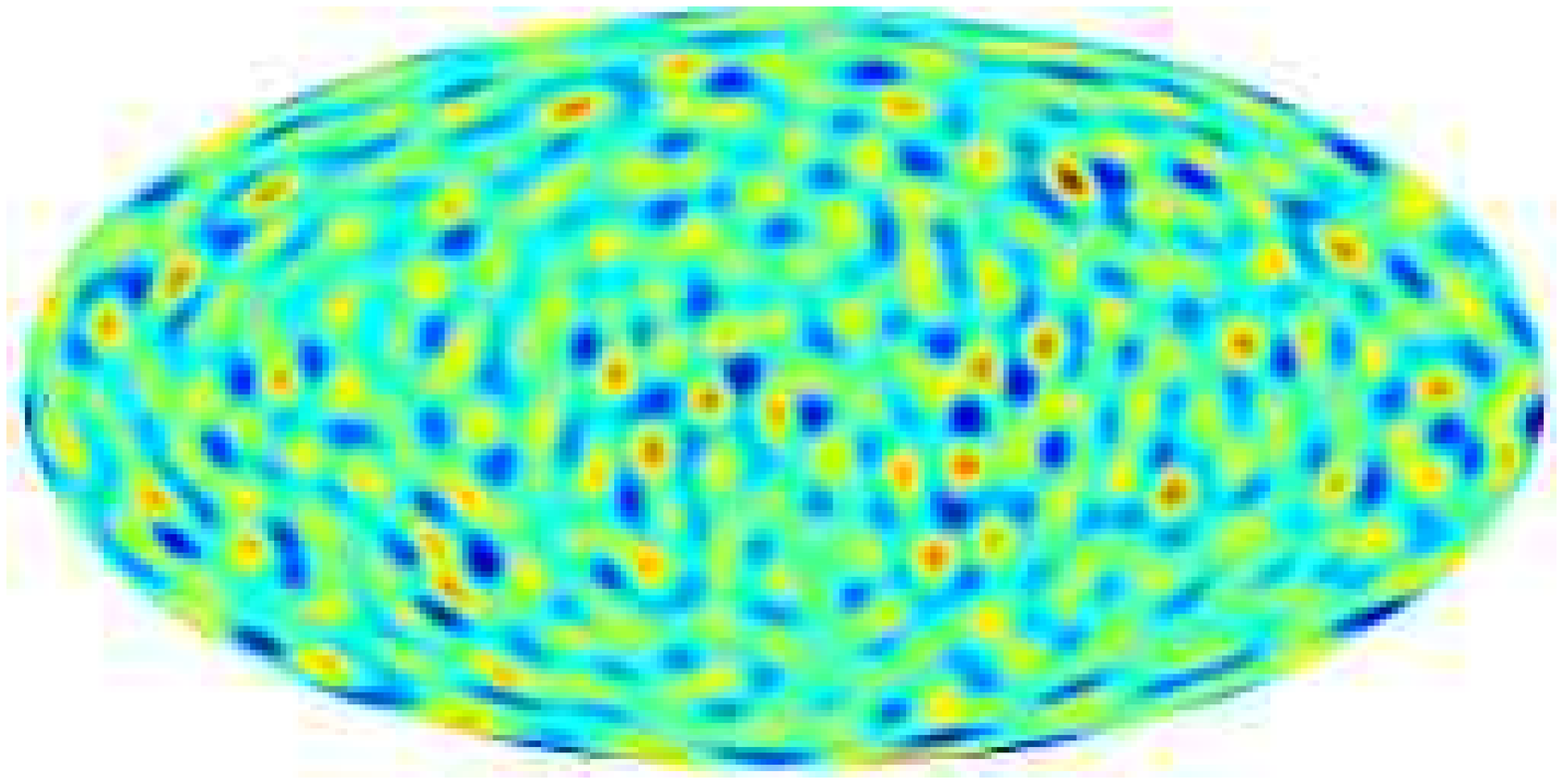}
    \includegraphics[width=\widthsky]{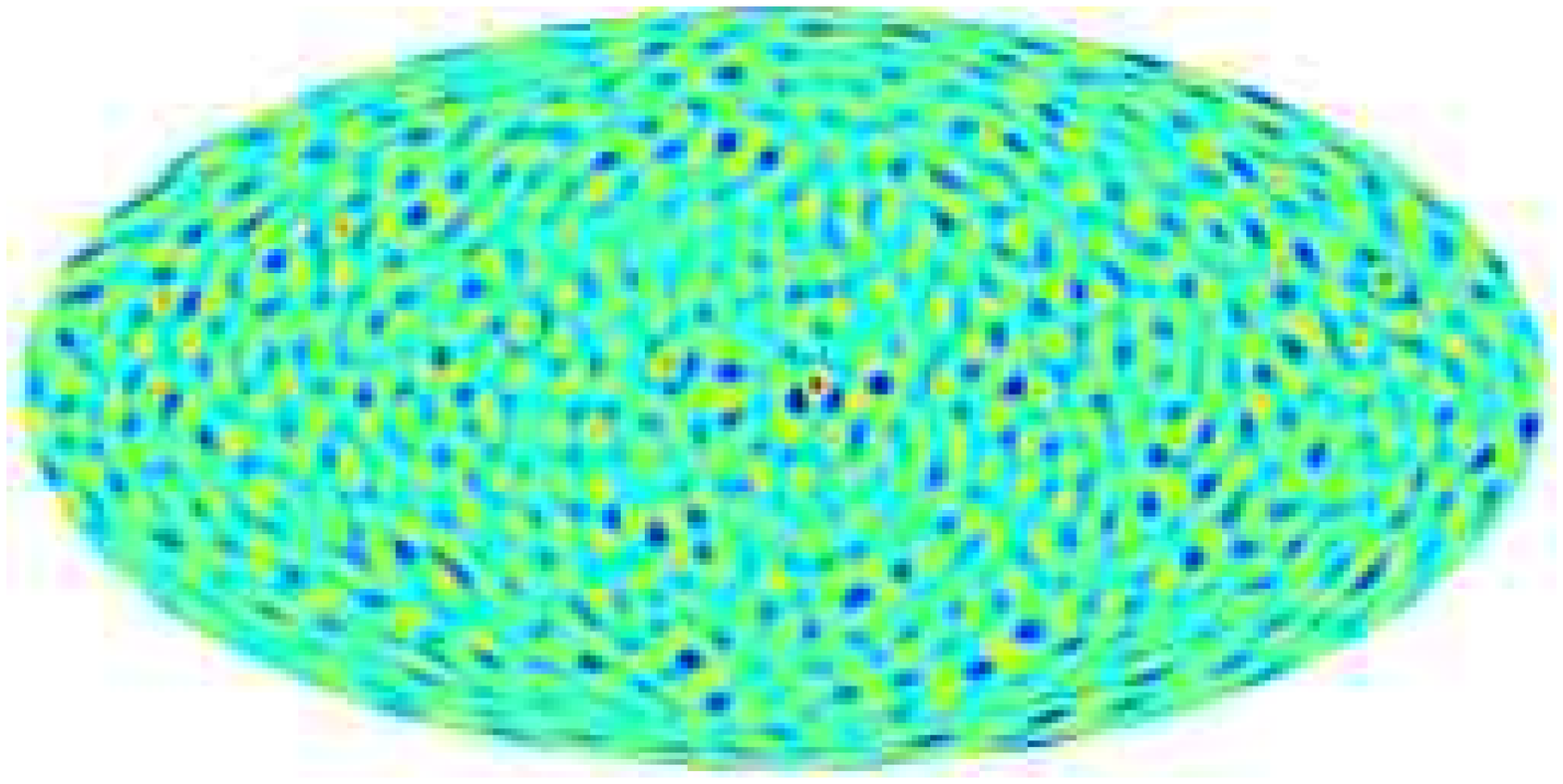}
    \includegraphics[width=\widthsky]{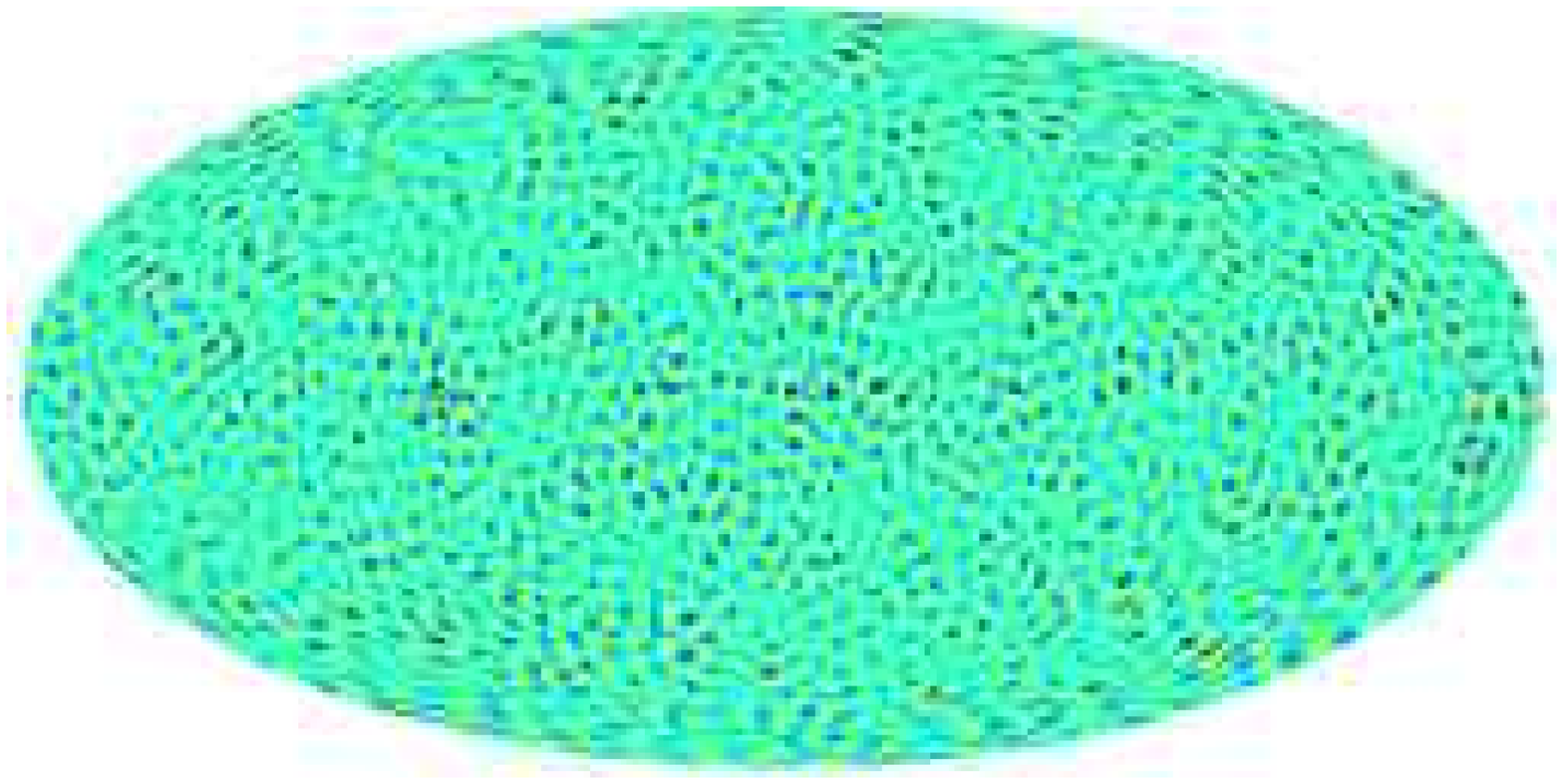}
    \includegraphics[width=\widthsky]{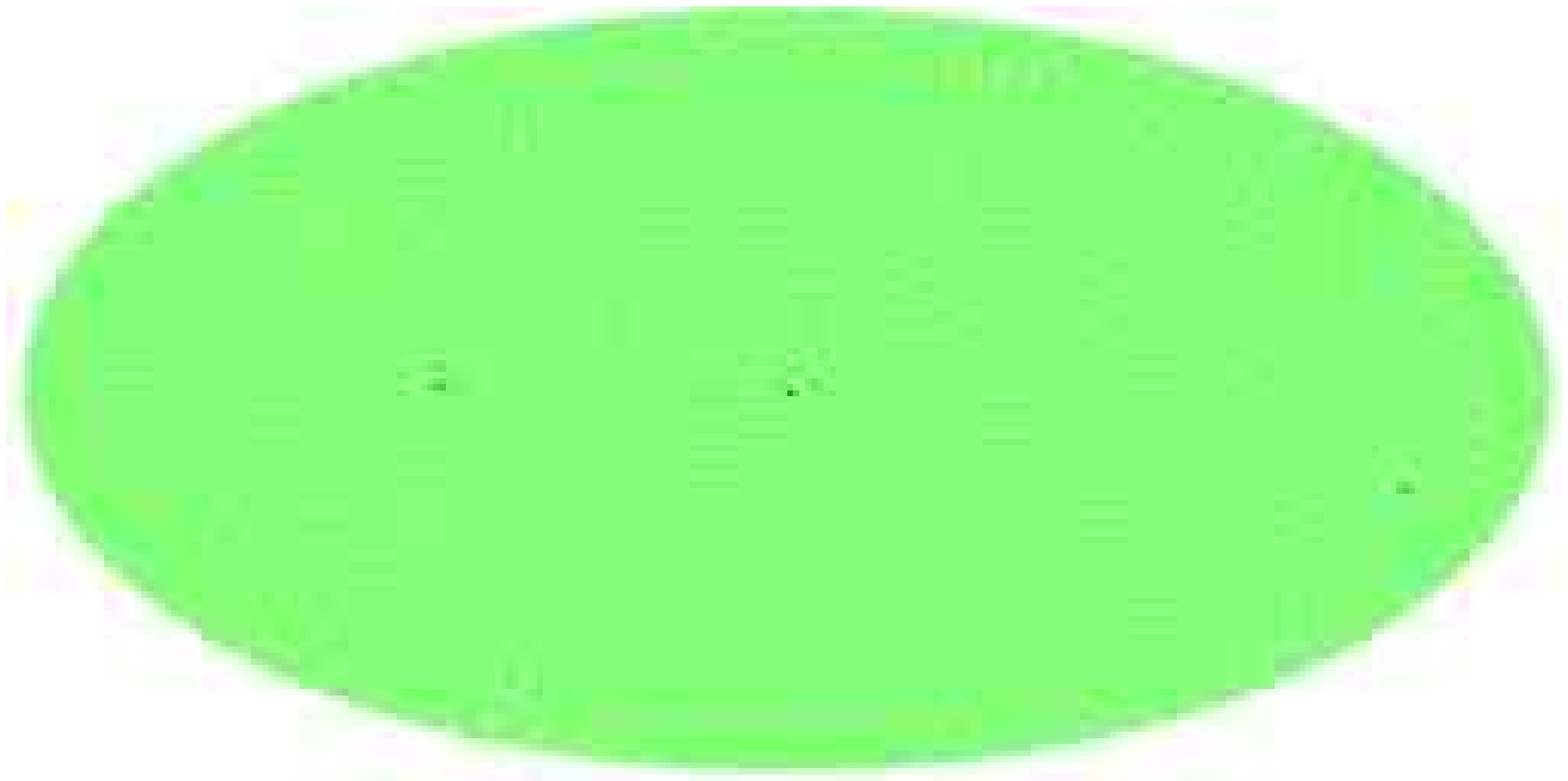}}
\caption{\label{fig:sph-liss} Input map of a
      CMB sky (from WMAP), and corresponding smoothed maps (with the
      spline filters of Figure \ref{fig:filt}).}
\end{figure}

\subsection{Needlet tight  frames}

Recall that a countable family of functions $\{f_n\}$ in a Hilbert
space $\mathcal H$  is a frame with frame bounds $C_1,C_2$ if
\begin{equation*}
  \forall g \in \mathcal H\; , \; C_1 \|g\|^2_{\mathcal H} \leq \sum_n | \langle
  g , f_n \rangle_{\mathcal H} |^2 \leq C_2 \|g\|^2_{\mathcal H} \; \; .
\end{equation*} It is a tight frame if we can choose $C_1 = C_2$.
Frames can be thought of  redundant ``bases'', and this redundancy can be exploited
for robustness issues. The tightness property is valuable in terms of numerical
stability \citep[see][Chap.3 and the references therein]{daubechies:1992}.

The construction that follows is from~\cite{narcowich:petrushev:ward:2006}. The
term \emph{needlet} was coined by~\cite{baldi:etal:2006a}.
Let $K$ be a finite index set and $\{\xi_k\}_{k \in K} \in \Sset^{|K|}$ a
set of quadrature points on the sphere, associated with a set
$\{\lambda_k\}_{k \in K} \in \Rset^{|K|}$ of quadrature weights.
\begin{defn}[Quadrature]
$\{(\xi_k,\lambda_k)\}_{k\in K}$ is said to provide an exact Gauss quadrature
formula at degree $\lmax$ if
\begin{equation*}
  \forall X\in\bigoplus_{\ell=0}^{\lmax}\Hset_\ell, \;
  \int_{\Sset} X(\xi) \dd \xi
  =
  \sum_{k\in K} \lambda_k X(\xi_k).
\end{equation*}
This quadrature formula is said positive-weight if $\lambda_k > 0, k \in K$.
\end{defn}
\begin{rem}
  We refer to \cite{Doroshkevich+2005} for an example of a proper choice of
  quadrature points and weights that fulfils this property (called
  GLESP). Other pixelization schemes such as HEALPix \citep{Gorski+2005} fulfil
  approximately this property with a number of points of order $C\lmax^2$ and
  quadrature weights of order $\frac1{C\lmax^2}$ for some positive constant
  $C$.
\end{rem}

Suppose that the window functions $h^{(j)}$ are non-negative and with
finite spectral  support. Define
\begin{equation}
\label{eq:bsqrth}
\forall \ell \in \Nset, \; b^{(j)}_\ell := \sqrt{h^{(j)}_\ell}   
\end{equation}
and $d^{(j)} := \max\{\ell:h^{(j)}_\ell\neq0\}$ (in
the $B$-adic case, $d^{(j)}=B^{j+1}$).
For each scale $j$, we have a pixellization
$\{\xi_k^{(j)},\lambda_k^{(j)}\}_{k \in K^{(j)}}$. 
\begin{defn}[Needlets and Needlet coefficients]
  \label{def:needlet}
  For every $j \in \mcJ$ and every index
  $k\in K^{(j)}$ the function 
\begin{equation}
  \label{eq:defneedlet}
  \psi_k^{(j)}(\xi)=\sqrt{\lambda_{k}^{(j)}}\sum_{\ell=0}^{d^{(j)}}
  b^{(j)}_\ell L_\ell (\xi \cdot \xi_k^{(j)} ),
\end{equation}
 is called a \emph{needlet}.
For $X \in \Hset$, the inner products $\langle
X,\psi_k^{(j)} \rangle$ are called \emph{needlet coefficients} and are
denoted $\beta_k^{(j)}$. 
\end{defn}
Up to a rotation of the sphere putting $\xi_k^{(j)}$ on the North pole and to
the multiplicative term $\sqrt{\lambda_k^{(j)}}$, all the needlets of a given scale
$j$ have exactly the
same shape.  In particular, they are axisymmetric.  When $\ell\mapsto
b^{(j)}_\ell$ is sufficiently smooth, one gets the intuition
from~(\ref{eq:defneedlet}) that the needlet $\psi_k^{(j)}$ is localized around
$\xi_k^{(j)}$.

The following Proposition state that the harmonic smoothing operation defined
by~(\ref{eq:defsmoothing}) can be seen as the decomposition of $\Hset$ on the
needlets family built with~(\ref{eq:bsqrth}), and that this family is a tight
frame.  It is a straightforward adaptation of
\citet[][Proposition~2.3]{baldi:etal:2006a}. 

\begin{prop} \label{prop:tightframe} 
Let $j\in\mcJ$. Assume that 
  $\{(\xi_k^{(j)},\lambda_k^{(j)})\}_{k \in K^{(j)}}$ provides an exact and
  positive-weight quadrature formula at degree $2 d^{(j)}$. Then
$$
\Psi^{(j)} X = \sum_{k\in K^{(j)}} \beta_k^{(j)} \psi_k^{(j)}.
$$
Assume that for any $j \in
  \mcJ$, $\{(\xi_k^{(j)},\lambda_k^{(j)})\}_{k \in K^{(j)}}$ provides an exact
  and positive-weight quadrature formula at degree $2 d^{(j)}$. Under the exact
  reconstruction condition~(\ref{reconstruction}),
 $$\forall X \in \Hset , \; X = \sum_{j\in\mcJ}\sum_{k \in K^{(j)}} \beta_k^{(j)} \psi_k^{(j)} \; \text{ and
 }\; \|X\|^2 = \sum_{j\in\mcJ}\sum_{k \in K^{(j)}} | \beta_k^{(j)}
 |^2 \; .$$
\end{prop}

\paragraph*{Remark on Terminology}
The analysis of an input field $X$ in the way described above is called
\emph{filtering}. This filtering has two equivalent expressions, in the spatial and in
the spectral domains; see the convolution formula~(\ref{eq:convolution}). These
expressions involves two ``dual'' mathematical objects : the functions
$h^{(j)}$ and $b^{(j)}$ of the frequency $\ell$, called \emph{window functions}
or \emph{spectral windows}, and the spherical functions $\psi_k^{(j)}$ called
\emph{needlets}, which are nothing else but the rotated axisymmetric functions
built from the Legendre transform of $b^{(j)}$ (see Definition~\ref{def:needlet}). We
call \emph{filter} either of the two above objects, when the domain (spatial or
spectral) is not specified.

\subsection{Generalized needlet frames}

We are concerned with the development of a flexible spectral analysis on the
sphere which remains practical at high resolution. The forecoming CMB
experiment \emph{Planck}\footnote{see \texttt{www.rssd.esa.int/Planck/}.} will provide
50 mega-pixel maps with accuracy such that multipole moments will be reliable up to
$\ell\simeq4000$.

For maximum flexibility, we shall consider constructions which are not
necessarily dyadic nor $B$-adic. This is motivated by  applications, as
described in the Introduction. Moreover, we will design analysis frames which
will not be necessarily tight. Their dual frames will be the corresponding
reconstruction frames. This allows fine tuning of the localization properties
of the decomposition functions but it is also well known that it does not
ensure similar properties for the reconstruction functions.  Nevertheless, for
the application goals discussed in the introduction, we will design strictly
band-limited needlets with support $L^{(j)} :=
[\ell_{\min}^{(j)},\ell_{\max}^{(j)}]$, $\ell_{\min}^{(j)}>0$ if $j \geq 0$.
Then the subsequent ``wavelet design'' operations will be performed in the
harmonic domain.

Since the needlet coefficients $\beta_k^{(j)}$ and $\beta_{k'}^{(j')}$ of a
Gaussian stationary (\ie{} isotropic) field are independent if $L^{(j)}\cap
L^{(j')} = \emptyset$, the bands $L^{(j)}$ are chosen to overlap as little as
possible. Other choices are possible; for instance \cite{Starck+2006} take
overlapping spectral windows supported on $[0, 2^j]$.

The three ingredients for our spherical ``multi-resolution'' approach are
harmonic-space implementation, dual wavelet frames and spectral window design.
In this subsection, we briefly describe the first two elements. In
Section~\ref{sec:crit-filt-design}, we go into the theory and practice
of window design.

\subsubsection{Dual  frames}

Proposition~\ref{prop:tightframe} shows that the needlets of
Definition~\ref{def:needlet} with~(\ref{eq:bsqrth}) can be used in both
analysis (or decomposition) and synthesis (or reconstruction). This accounts to
say that the needlet frame is its own dual frame.  We choose to keep the
Definition~\ref{def:needlet} of the needlets and associated coefficients but to
relax condition~(\ref{eq:bsqrth}). By sacrificing the tightness of the frame,
we gain much freedom in the design of the spectral windows. Also, the precise
space-frequency picture provided by the needlet construction is preserved.

From any windows family $(b^{(j)})_{j \in \mcJ}$ such that $\forall \ell \in \Nset, \sum_{j
  \in \mcJ}\left(b^{(j)}_\ell\right)^2 > 0$, define the synthesis windows
$\tilde{b}^{(j)}$ by
\begin{equation}
  \label{eq:defbtilde} 
  \forall j \in \mcJ, \;  \forall \ell \in \Nset, \; \tilde{b}^{(j)}_\ell
  = \frac{b^{(j)}_\ell}{\sum_{j' \in \mcJ}\left(b^{(j')}_\ell\right)^2} \,
\end{equation}
and put $h^{(j)} := \tilde b^{(j)} b^{(j)}$ so that (\ref{reconstruction})
easily follows.  We retain Definition~\ref{def:needlet} for the
\emph{decomposition} needlets and needlets coefficients and further define the
\emph{reconstruction} needlets as
\begin{equation}
   \label{eq:defreconstructionneedlet}
  \tilde \psi_k^{(j)}(\xi)=\sqrt{\lambda_{k}^{(j)}}\sum_{\ell=0}^{d^{(j)}}
  \tilde b^{(j)}_\ell L_\ell (\xi \cdot \xi_k^{(j)} ) \ .
\end{equation}
\begin{prop}\label{prop:frame}
  Assume that there exists positive constants $ C_1, C_2$ such that
  \begin{equation}
    \label{eq:boundssumblj}
    \forall \ell \in \Nset ,  C_1 \leq \sum_{j \in \mcJ}|b_\ell^{(j)}|^2 \leq
    C_2 \ .
  \end{equation}
  Assume that for any $j \in \mcJ$, the set $\{(\xi_k^{(j)},\lambda_k^{(j)})\}_{k
    \in K^{(j)}}$ provides an exact and positive-weight quadrature formula at
  degree $2 d^{(j)}$.  Then, under the exact reconstruction
  condition~(\ref{reconstruction}), the family $\{\psi_k^{(j)}\}$ is a frame
  with frame bounds constant $C_1$ and $C_2$. Its dual frame is the family
  $\{\tilde \psi_k^{(j)}\}$. In particular
  \begin{equation}
    \label{eq:dualframes}
\forall X \in \Hset , \; X \stackrel{\Hset}{=} \sum_{j\in\mcJ}\sum_{k \in K^{(j)}} \beta_k^{(j)} \tilde\psi_k^{(j)} \; \text{ and
 }\; \|X\|^2 = \sum_{j\in\mcJ}\sum_{k \in K^{(j)}} \tilde\beta_k^{(j)}\beta_k^{(j)} \; ,
  \end{equation} 
with $\tilde\beta_k^{(j)}:=\langle X,\tilde\psi_k^{(j)}\rangle$.
\end{prop}

Define the analysis, synthesis and smoothing operators at
scale $j \in \mcJ$ by $\Phi^{(j)} = \sum_\ell b^{(j)}_\ell \Pi_\ell$, $\tilde \Phi^{(j)}
= \sum\limits_\ell \tilde b^{(j)}_\ell \Pi_\ell$ and $\Psi^{(j)}
=\tilde\Phi^{(j)}\Phi^{(j)}$, respectively. Then, the exact
reconstruction formula $\sum \Psi^{(j)} = \textrm{Id}$ holds true.

\begin{figure}[htbp]
  \centering  \includegraphics{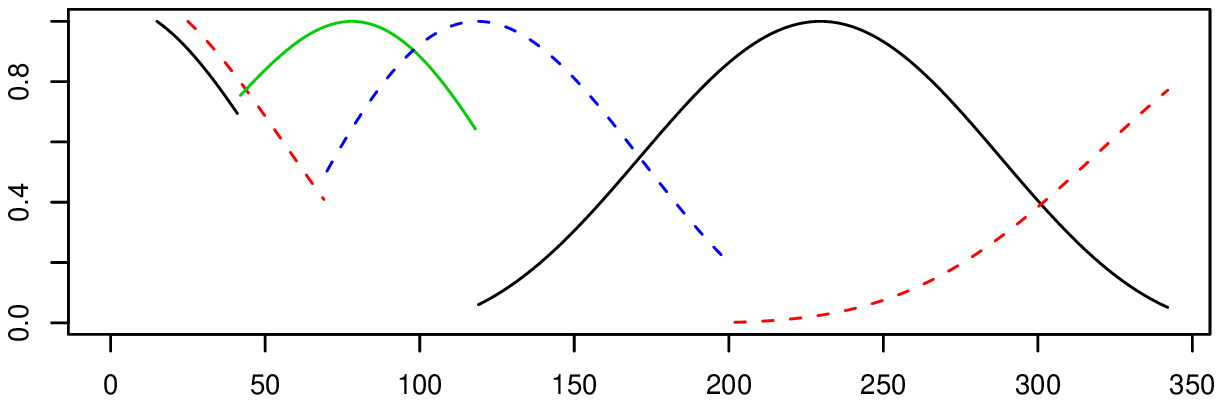}
  \includegraphics{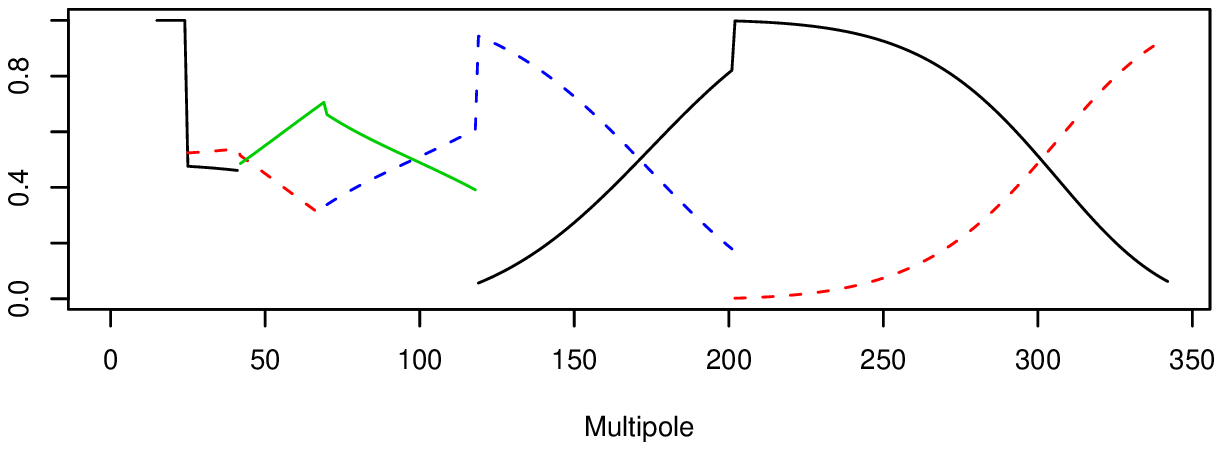}
  \caption{ \label{fig:dualframes} $B$-adic analysis (top) and
    corresponding synthesis (bottom) window functions ($j=6,\dots,11\
    ;\ B=1.7$).}
\end{figure}

An example of an analysis/synthesis windows family following this
scheme is displayed in Figure~\ref{fig:dualframes}, in which we took optimally
concentrated PSWF (see Section~\ref{sec:crit-filt-design}) functions for
analysis.  It illustrates the fact that this choice does not lead to well
localized synthesis needlets (as their spectral shapes are non
smooth). However, this may not be a shortcoming if one is interested in the
needlet coefficients $\beta_k^{(j)}=\langle X,\psi_k^{(j)}\rangle$ \emph{per
  se}, which reflect the local properties of the field $X$.

\subsubsection{Practical computation of  needlet coefficients}

Evaluation of inner products $\langle X,\psi^{(j)}_k \rangle $ in the direct
space is practically unfeasible from a pixelized sphere at high
resolutions. The needlet coefficients $\beta^{(j)}_k$ are thus computed via
direct and inverse harmonic transforms as a consequence of the
following Proposition.

\begin{prop} \label{prop:expneedcoeff}
  The needlet coefficients  verify
  $\beta^{(j)}_k = \sqrt{\lambda^{(j)}_{k}}\Phi^{(j)} X(\xi^{(j)}_k).$
\end{prop}

The computation of the smoothed field $\Phi^{(j)} X$ is performed in
the harmonic domain by multiplying the multipole coefficients $a\lm$
of $X$ by the factors $b^{(j)}_\ell$. Finally, the needlet
coefficients $\beta^{(j)}_k$ are retrieved as the values of
$\Phi^{(j)} X$ at the points $\xi^{(j)}_k$ up to a multiplicative
term. Starting from the field $X$ sampled at some quadrature points,
this operation is summed up by the diagram
\begin{equation}
  \label{eq:diagram}
  \{X(\xi_k)\}_{k \in K}
  \stackrel{\textrm{SHT}}{\longrightarrow}\{a\lm\}\lm
  \stackrel{\times}{\longrightarrow} \{b^{(j)}_\ell a\lm\}\lm
  \stackrel{\textrm{SHT}^{-1}}{\longrightarrow}\left\{(\lambda^{(j)}_k)^{-1/2}\beta^{(j)}_k\right\}_{k
    \in K^{(j)}}  
\end{equation}
whereas the synthesis operation is summed up by
$$\left\{(\lambda^{(j)}_k)^{-1/2}\beta^{(j)}_k\right\}_{k \in K^{(j)}}
\stackrel{\textrm{SHT}}{\longrightarrow} \{b^{(j)}_\ell a\lm\}\lm
\stackrel{\times}{\longrightarrow} \{\tilde b^{(j)}_\ell b^{(j)}_\ell a\lm\}\lm
\stackrel{\textrm{SHT}^{-1}}{\longrightarrow}\{\Psi^{(j)} X(\xi^{(j)}_k)\}_{k\in
  K^{(j)}}
$$

Standard pixelization packages, such as HEALPix, GLESP or
SHTOOLS\footnote{available at
  http://www.ipgp.jussieu.fr/$\sim$wieczor/SHTOOLS/SHTOOLS.html} come with
optimized implementations of the direct and inverse Spherical Harmonic
Transforms. For example, in the HEALPix scheme, pixels are located on rings of
constant latitude, allowing for fast SHT.
This makes the computation easy and tractable even at high
resolution. The needlet coefficients at a given scale $j$ can be
visualized as a pixelized map. If the quadrature weights
$\{\lambda_k^{(j)}\}$ are equal, the smoothed maps of
Fig.~\ref{fig:sph-liss}, which are the outputs of the processing
(\ref{eq:diagram}), provide a precise and easily interpretable picture
of the space-frequency analysis.

\begin{rem}
  The quadrature points and weights $\{(\xi^{(j)}_k,\lambda^{(j)}_k)\}_{k \in K^{(j)}}$ use to
  define the needlet coefficients $\beta^{(j)}_k$ and to sample the smoothed
  field $\Psi^{(j)} X$ may be chosen identical to $\{(\xi_k,\lambda_k)\}_{k \in K}$ used
  to sample the input field $X$. However, for data compression and
  computational efficiency, one can consider alternatively to take
  the minimal $K^{(j)}$ providing an exact positive-weight quadrature formula at a proper
  degree.
\end{rem}

\section{Design of optimally localized wavelets}
\label{sec:crit-filt-design}

In this section, we define some criteria to compare the window
profiles.  Some of them are easily optimized, others are only
investigated numerically.  We first give some examples of generic
needlet profiles we can think of (Section~\ref{sec:examples}). Then,
we restrict ourselves to a single scale $j$ and an associated band
$L:=[\lmin,\lmax]$. The superscript $(j)$ will be omitted in the
notations when no confusion is possible. We present the $\Lset^2$
(Section~\ref{sec:lset2}) and statistical
(Section~\ref{sec:statistical-criteria}) criteria, with practical
implementation details on their optimizations.

\subsection{Examples.}
\label{sec:examples}

\cite{narcowich:petrushev:ward:2006} have derived the following theoretical
bound that controls the decay of the needlets.  In the $B$-adic case, if the
function $\mathsf b := \sqrt{\mathsf{h}}$ defining the analysis spectral window
is $M$-times continuously differentiable,
\begin{equation*}
  |\psi_k^{(j)}(\xi)|\leq \frac{C\
    B^{j-1}}{1+\bigl(B^{j-1}\arccos(\xi \cdot \xi_k^{(j)})\bigr)^M}
\end{equation*}
for some constant $C=C(\mathsf b)$.  This condition still allows a wide range
of possibilities for designing the function $\mathsf b$. Without restricting
ourselves to the $B$-adic case, we implemented solutions to optimize in
practice, non asymptotically, the shape of windows $b^{(j)}$ regarding some
applications.

To illustrate the kind of aspects we are concerned with, we compare in
Figure~\ref{fig:examples} the azimuthal profiles (in the spatial
domain) of various axisymmetric needlets. The needlets are built from
window functions $b^{(j)}$ via relation~(\ref{eq:defneedlet}) and
$\xi_k = (0,0)$, \ie{} they are centered on the North pole, and then
are considered as functions of $\theta$ only.  This illustration is
restricted to the $9^{\textrm{th}}$ dyadic scale, \ie{} frequencies in
the band $L := [256,1024]$.  We shall compare heuristically five
families of window functions. Note that the last two are not limited
to band $L$.
\begin{enumerate}
\item Square roots of splines of various orders. For any odd integer
  $M$, there exists a spline function $\textsf h$ of order $M$,
  non-negative, compactly supported on $[\frac12,2]$ and such that the
  $h_\ell^{(j)}$'s defined by (\ref{eq:Badic})
  verify~(\ref{reconstruction}). It remains to define
  $b_\ell^{(j)}=\sqrt{h_\ell^{(j)}}$.
\item Best concentrated Slepian functions in caps of various radii (cf
  Section~\ref{sec:lset2}). The window function $b^{(j)}_\ell$ is the
  minimizer of the criterion~(\ref{eq:optimprob2}). It is band-limited
  on $L$ and optimally concentrated in a polar cap
  $\Omega_{\theta_0}=\{\xi:\theta\leq\theta_0\}$), $\theta_0$ being a
  free parameter.
\item Denote $G$ a primitive of the $C^\infty$ function $g:x\mapsto
  e^{-\frac1{1-x^2}}1_{(-1,1)}(x)$ and put
  \begin{equation}
    \label{eq:ital}
    \mathsf{b}(x)=G(-8x+3)-G(-4x+3)
  \end{equation}
  and $b^{(j)}_\ell=\mathsf{b}\left(\frac{\ell}{2^j}\right)$. This
  window function is used in \cite{pietrobon:balbi:marinucci:2006}.
\item From the $B$-spline function of order 3
  \begin{equation}
    \label{eq:B3}
    B_3(x)=\frac1{12}(|x-2|^3-4|x-1|^3+6|x|^3-4|x+1|^3+|x+2|^3),
  \end{equation}
  form $\mathsf{b}(x)=\frac32(B_3(2x) - B_3(x))$ and define
  $b^{(j)}_\ell=\mathsf{b}\left(\frac{\ell}{2^j}\right)$. This window
  function is used by \cite{Starck+2006}.
\item The Mexican hat wavelet on the sphere is the function the
  stereographic projection of which on the Euclidean plane is the
  usual Mexican hat wavelet. It has the following close expression
  depending on some positive scale parameter $R$
  \begin{equation}
    \label{eq:mexhat}
    \psi_R(\theta) \propto (1 - 2R^2\tan^2(\theta/2)) \exp\{ - 2R^2\tan^2(\theta/2)\}.
  \end{equation}
  This wavelet is  popular in the astrophysics community
  \citep[see \eg{}][]{gonzalez-nuevo:etal:2006}. We have chosen $R=6.10^{-3}$
  such that the spectral window is almost zero for $\ell > 1024$. 
\end{enumerate}

\newlength{\widthfig}
\setlength{\widthfig}{4.3cm}

\begin{figure}[htbp]
\centering

\subfigure[Splines of order resp. 7, 15, 31 and 43.]{
  \includegraphics[width=\widthfig]{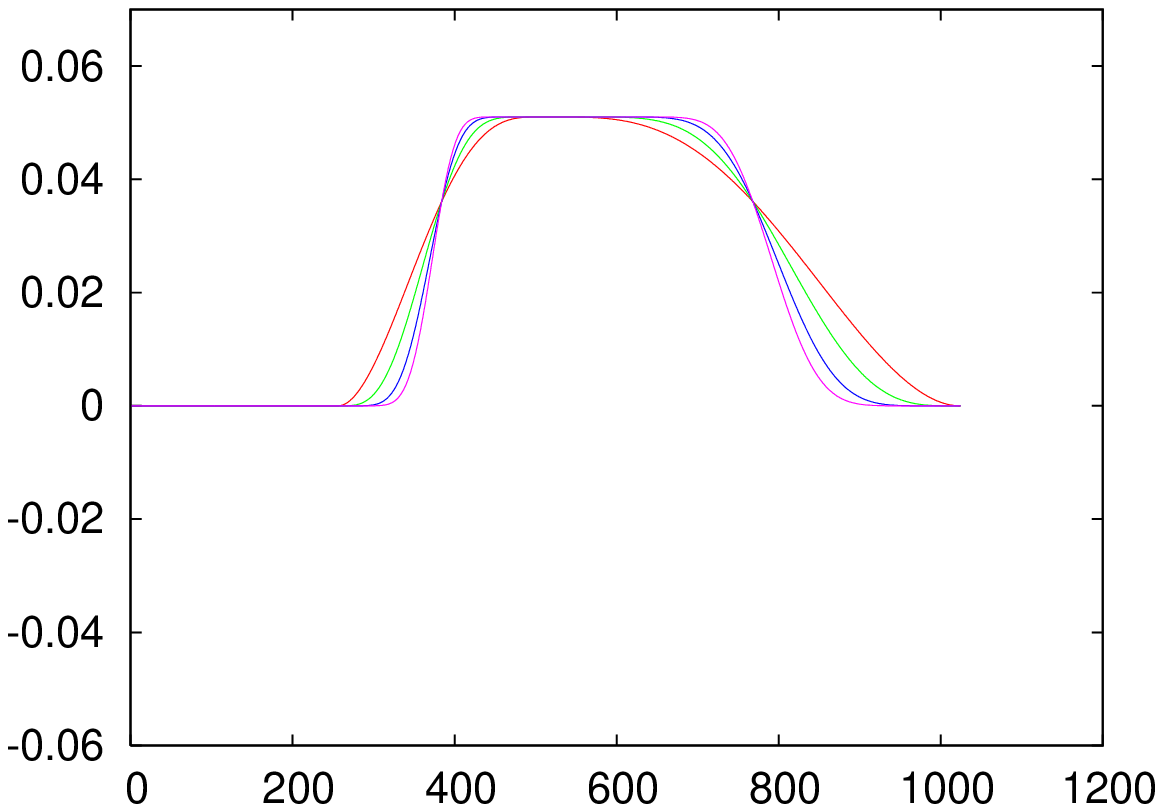}
  \includegraphics[width=\widthfig]{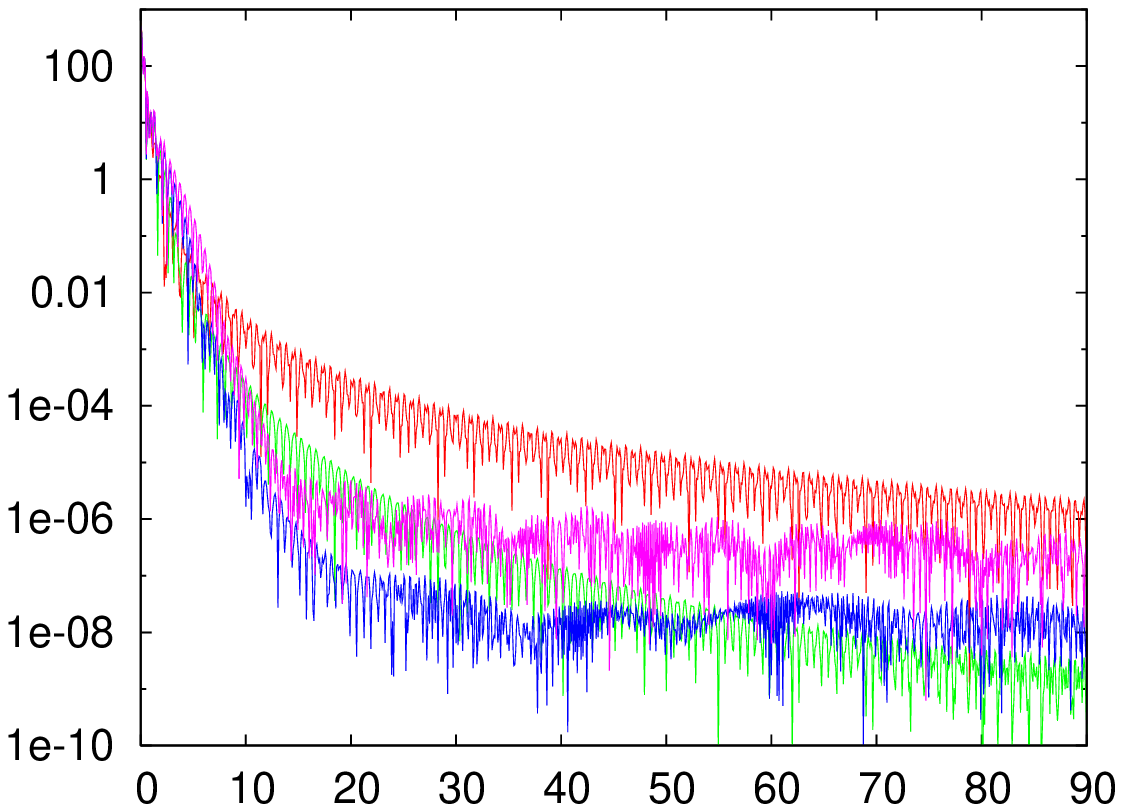}
  \includegraphics[width=\widthfig]{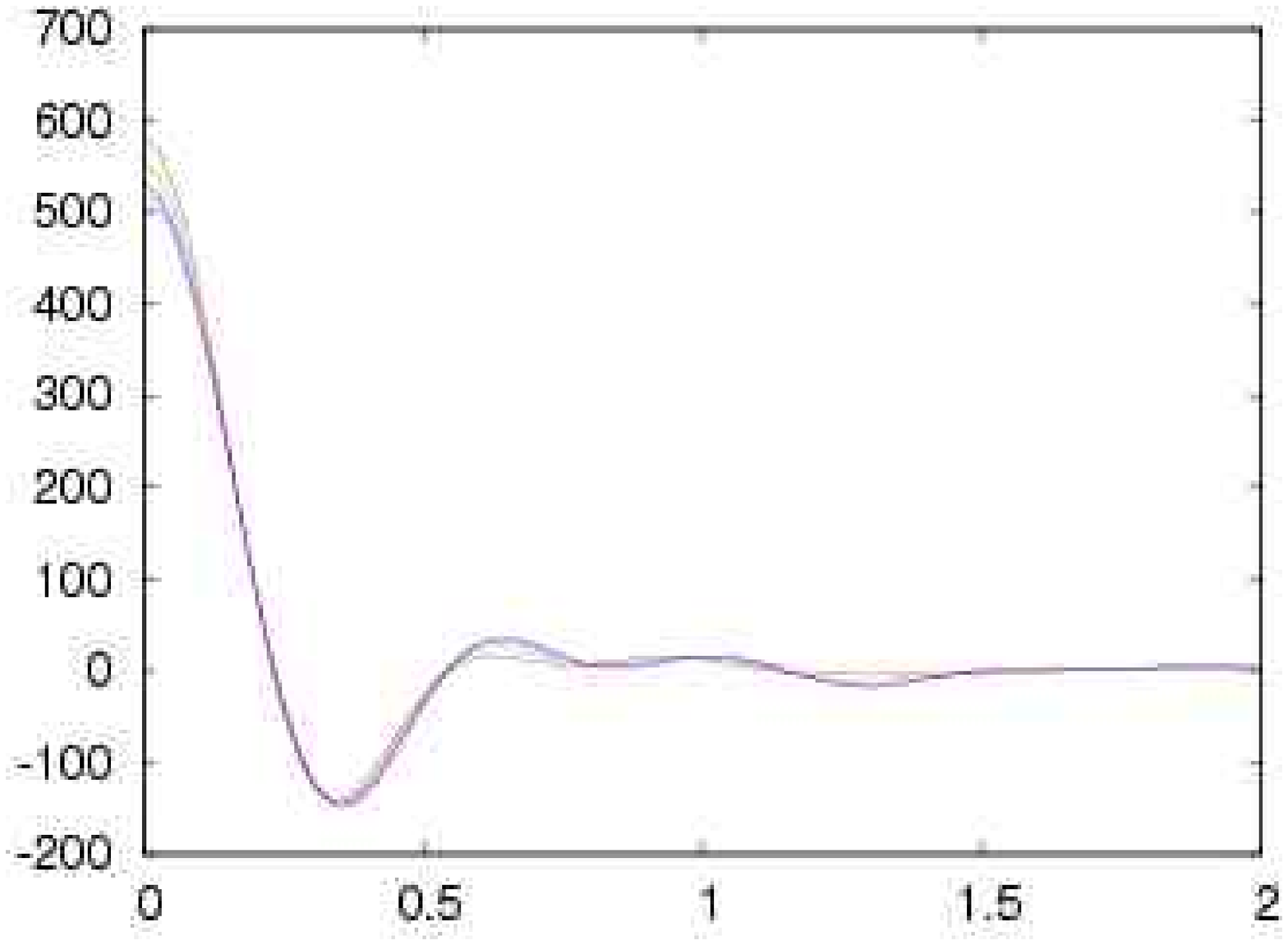}
}
\subfigure[PSWFs localized in polar caps of 0.5, 1, 1.5 and 5
degree opening]{
  \includegraphics[width=\widthfig]{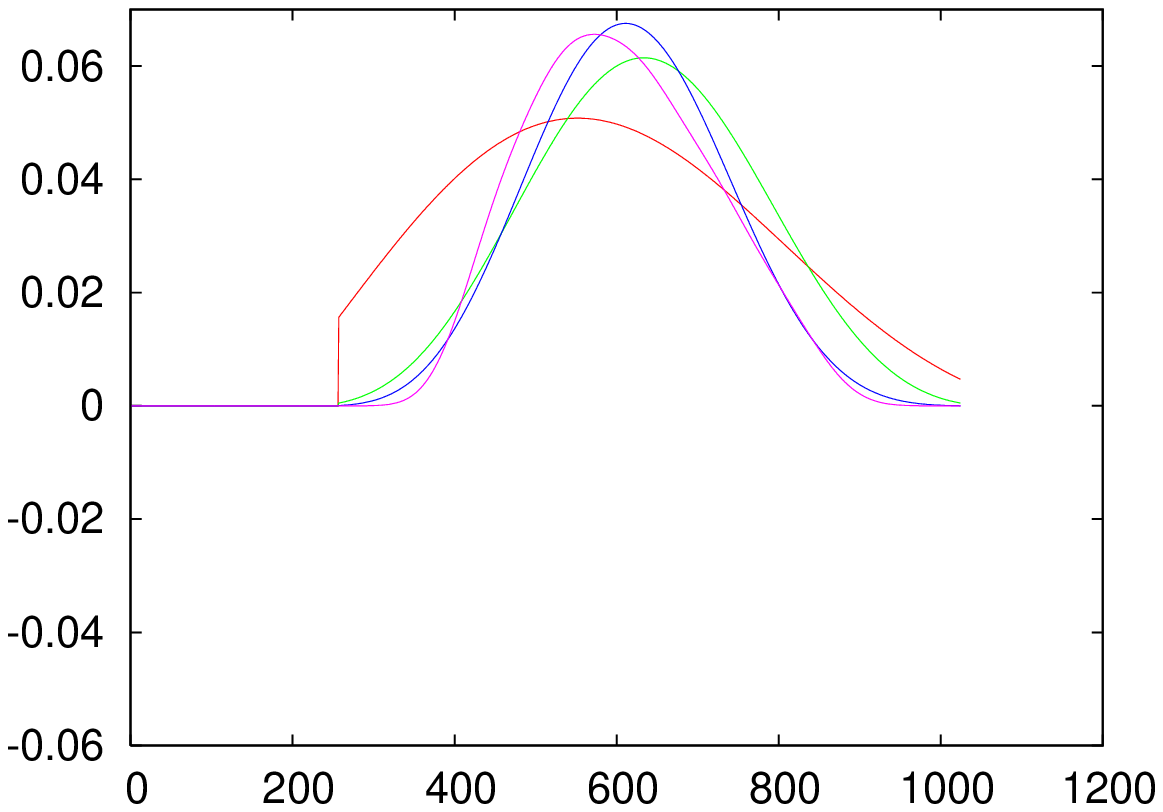}
  \includegraphics[width=\widthfig]{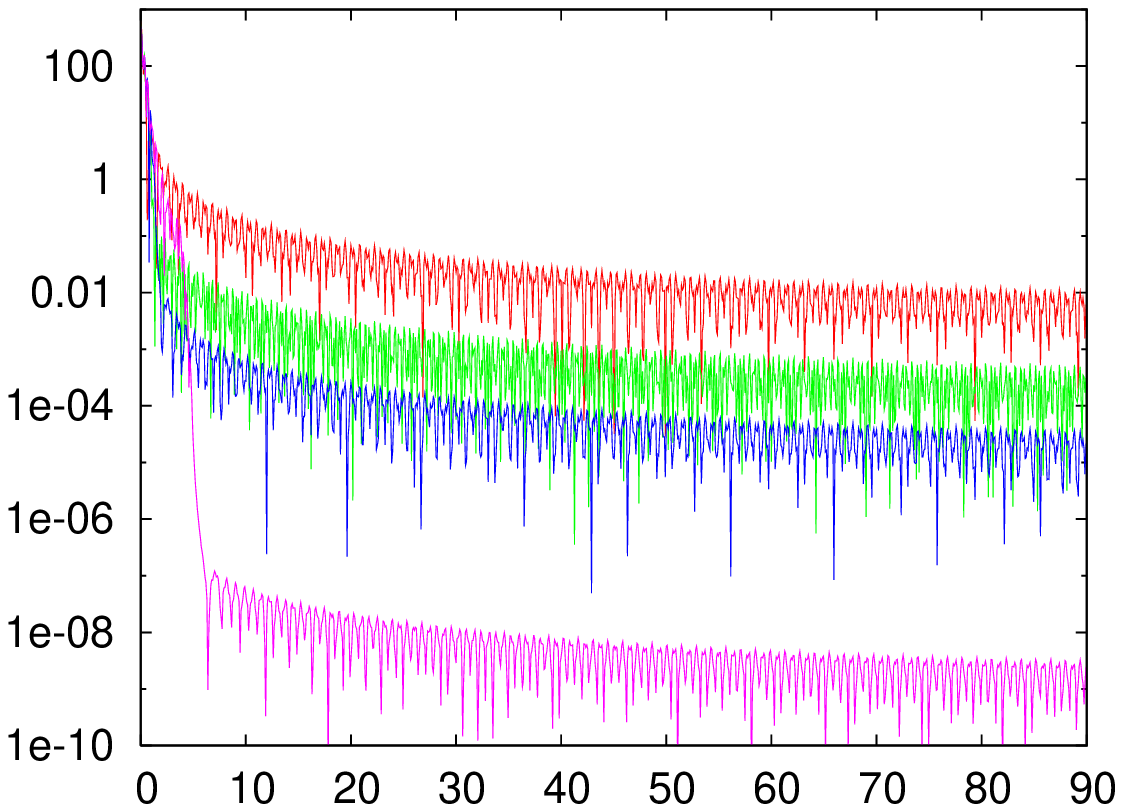}
  \includegraphics[width=\widthfig]{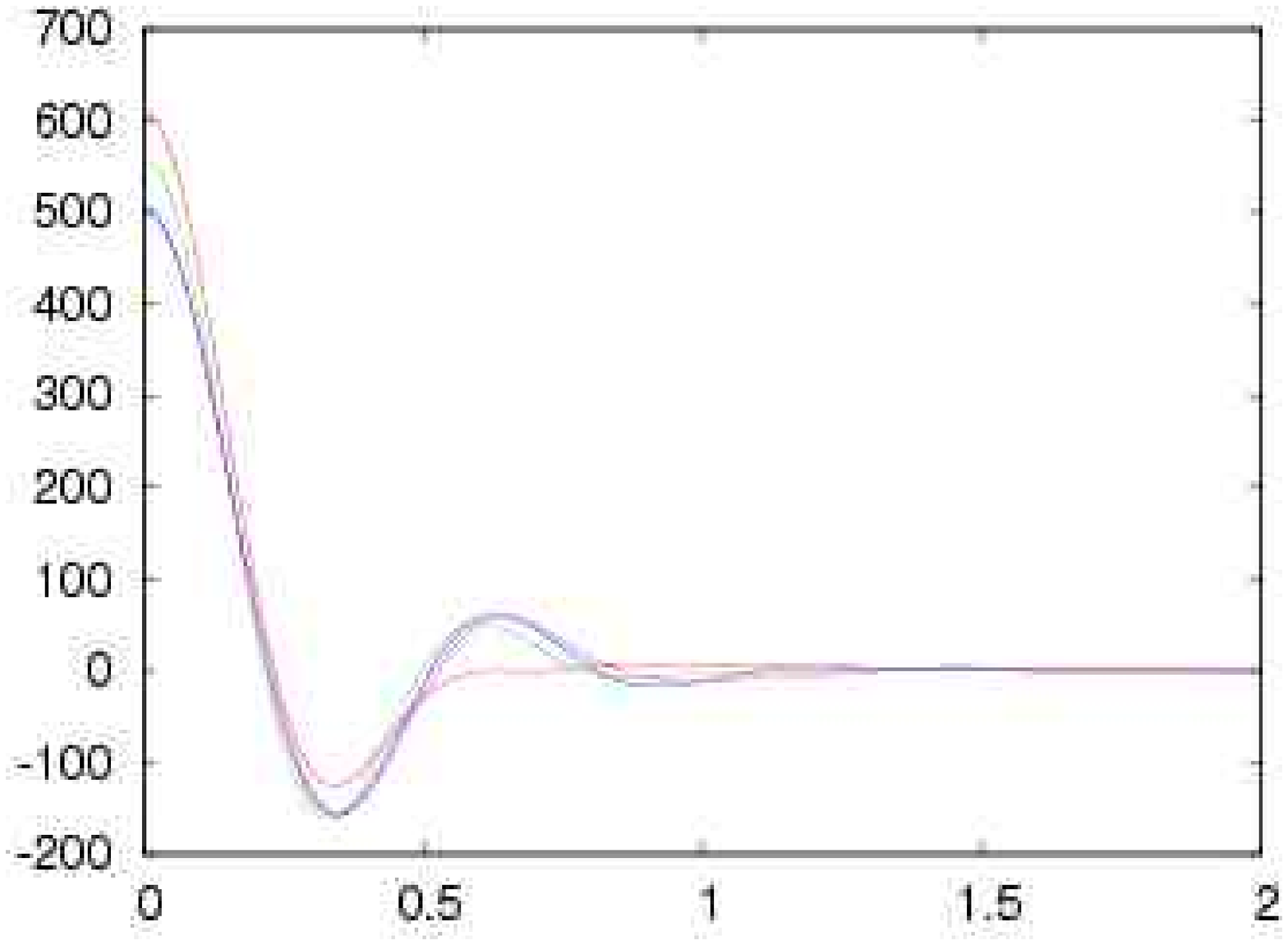}
}
\\
\subfigure[(red curve) Exponential function described in
Eq.~(\ref{eq:ital}), (green curve) $B$-spline function of
Eq.~(\ref{eq:B3}) and (blue curve) Mexican hat described in
Eq.~(\ref{eq:mexhat}).]{ 
  \includegraphics[width=\widthfig]{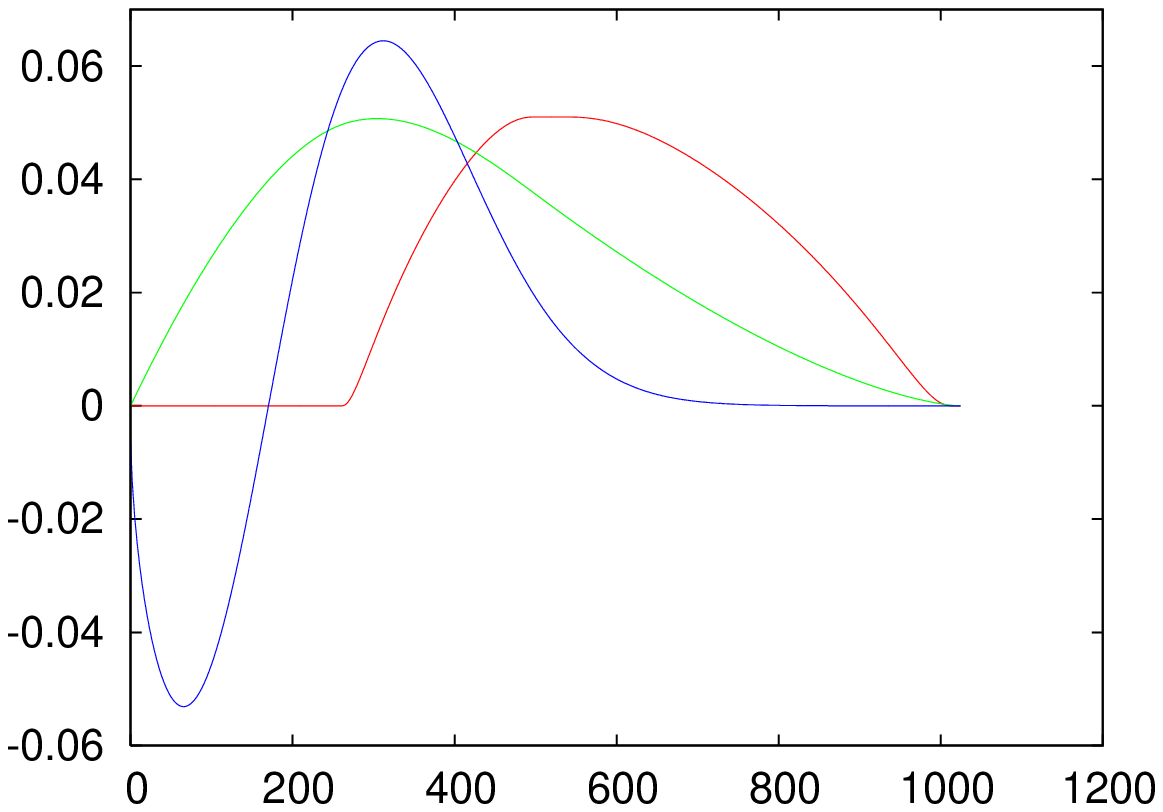}
  \includegraphics[width=\widthfig]{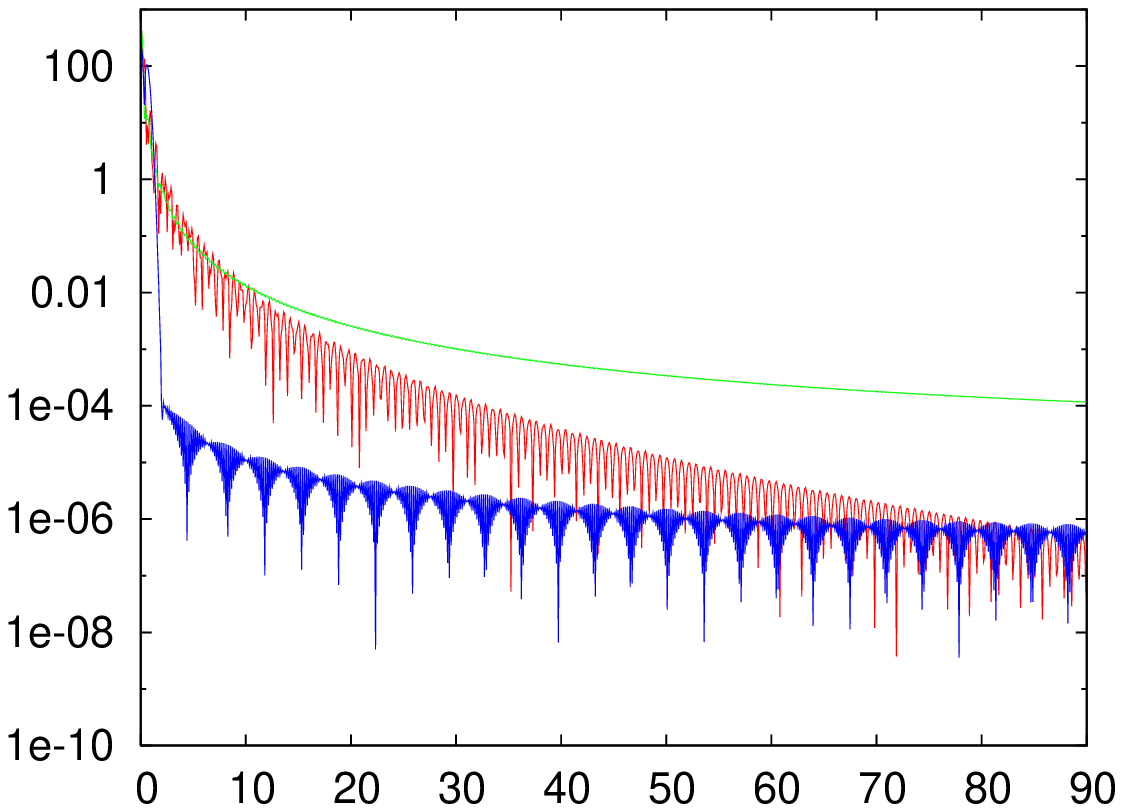}
  \includegraphics[width=\widthfig]{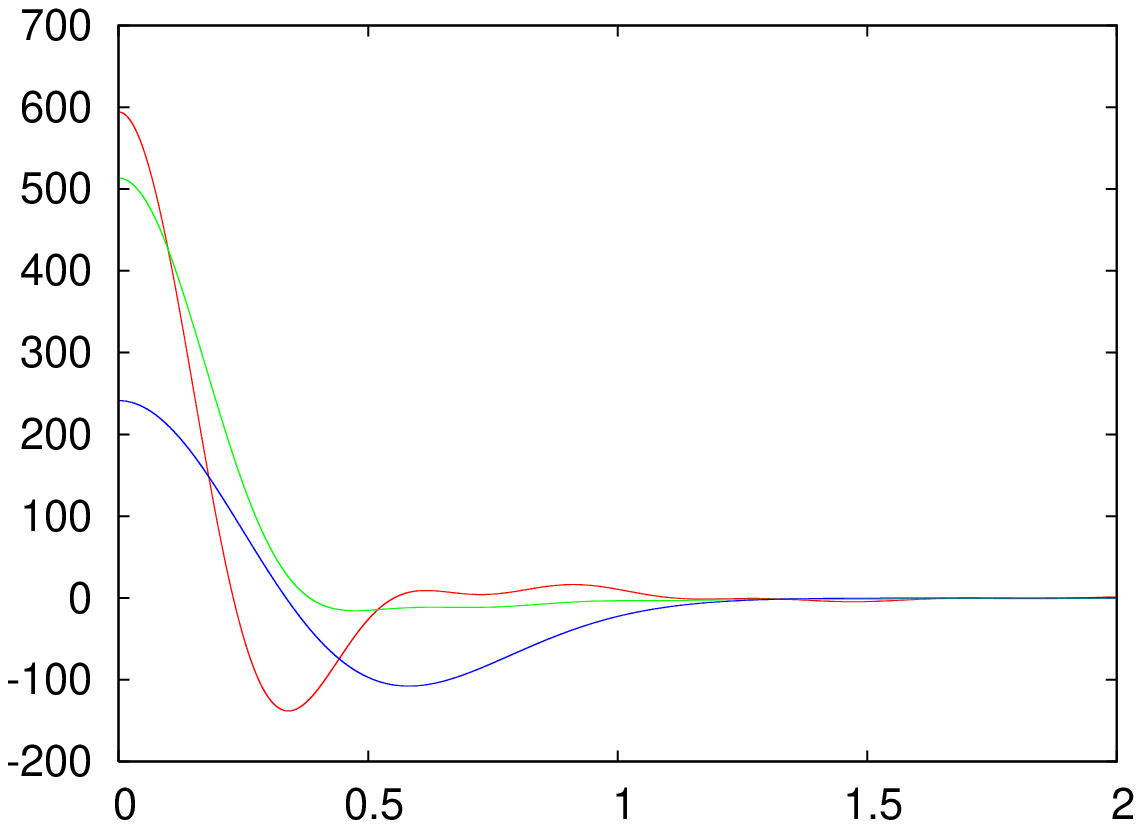}
}
\\
\caption{ \label{fig:examples} In left column, the shape of the
  spectral windows as a function of $\ell$. In middle and right
  columns, the profile of the filters is plotted in the spatial domain
  as a function of $\theta$ ($\theta$ in degrees) with logarithmic and
  linear scales respectively, to illustrate both the decrease of the
  tail of the needlets far from the North pole and the shape of their
  first bounces.}
\end{figure}

\subsection{$\Lset^2$-concentration  and variations}
\label{sec:lset2}

Our first attempt to achieve a good spatial localization of a needlet
is to optimize a $\Lset^2$-norm based criterion, adapting to the
sphere a problem that is well-known on the real line. In their seminal
work in the 1960s and 1970s, Slepian and his collaborators studied the
properties of prolate spheroidal wave functions (PSWFs) in the 1D case
of real functions
\citep[see][and the references therein]{Slepian83}.  PSWFs may be
defined as functions with optimal energy concentration in the time
domain, under some band-limitation constraint.  Equivalently, they are
the eigenfunctions of a time-frequency concentration kernel or the
solutions of a Sturm-Liouville differential equation. The
time-frequency concentration of PSWFs is understood in terms of
continuous Fourier transform on $\Rset$. A discrete version of this
theory, based on Fourier series coefficients, is derived in
\cite{Slepian+78}.

In the last few years, Walter and coauthors exploited these 1D PSWFs
to derive Slepian series (in \citealp{walter+2003}; see also
\citealp{Moore+2004}), and wavelets based on the best concentrated
PSWF \citep{walter+2004, Walter+2005}.

On the sphere, we shall only consider the equivalent of Discrete
PSWFs, following \cite{Simons+2006}. From a window function
$\{b_\ell\}$ with support $L$, define the axisymmetric
function $\psi$ by
\begin{equation}
  \label{eq:B_L}
  \psi(\xi) = \sum_{\ell \in L} b_\ell L_\ell(\cos\theta).
\end{equation}
The set of functions $\psi$ of the form (\ref{eq:B_L}) is denoted
$B_L\subset\Hset$.  
Given a spherical domain $\Omega$, consider the minimization, among
non-zero functions in $B_L$, of the criterion
\begin{equation}
  \label{eq:optimprob1}
  \mathcal{C}_\Omega(\psi)
  =
  \frac{\int_{\Sset\setminus\Omega} \psi^2(\xi)\dd\xi}{\int_{\Sset} \psi^2(\xi) \dd\xi}
  =
  1-\frac{\int_\Omega \psi^2(\xi) \dd \xi}{\int_{\Sset} \psi^2(\xi) \dd \xi}.
\end{equation}
This extension to the sphere of Slepian's concentration problem is studied in
details by \cite{Simons+2006} in the case $\lmin=0$.  We call PSWF (by abuse of
language) and denote $\psi^\star_{\Omega}$ a normalized minimizer for $
\mathcal{C}_\Omega(\psi)$.

\begin{figure}[htbp]
  \centering
  \includegraphics{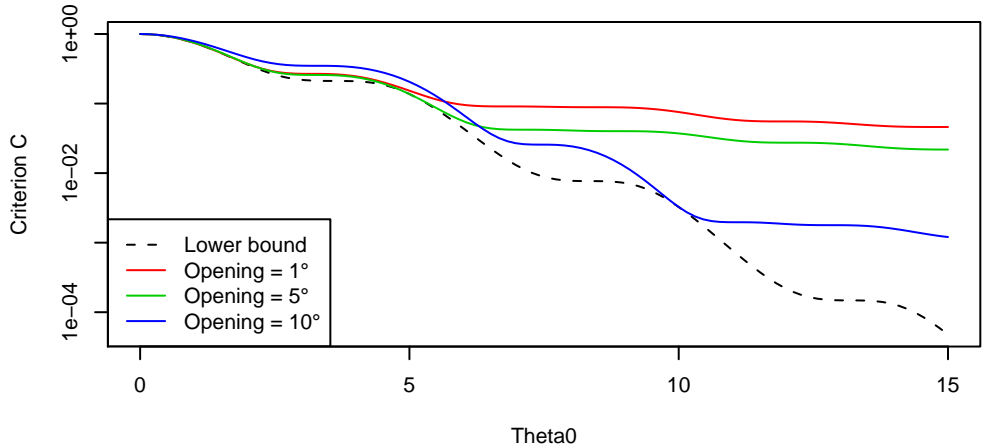}
  \caption{\label{fig:prol_local} Localization for $\Lset^2$-energy criterion
    of PSWFs, band-limited into $L=[33,64]$. The dashed line is the minimum
    of the criterion $\mathcal{C}_{\Omega_{\theta_0}}$ as a function of
    $\theta_0$ and the other ones are the values
    of $\mathcal{C}_{\Omega_{\theta_0}}(\psi)$ evaluated at  $\psi = \psi^\star_{\Omega_1}$, $\psi^\star_{\Omega_5}$ and
    $\psi^\star_{\Omega_{10}}$. }
\end{figure}

The criterion~(\ref{eq:optimprob1}) has a simplified expression when
$\Omega$ is axisymmetric.  Consider the polar cap
$\Omega_{\theta_0}=\{\xi:\theta\leq\theta_0\}$ and define the coupling
matrix $\bD = (D_{\ell,\ell'})_{\ell,\ell'\in L}$ by $$D_{\ell,\ell'}
= \frac{8\pi^2}{\sqrt{(2\ell+1)(2\ell'+1)}}\int_{\cos\theta_0}^1
L_\ell(z) L_{\ell'}(z) \dd z \ , $$and $$ \bar{\bb}(\psi) =
(\sqrt{\frac{2\lmin+1}{8\pi^2}}b_{\lmin},\dots,
\sqrt{\frac{2\lmax+1}{8\pi^2}}b_{\lmax}).$$ Then 
\begin{equation}
  \label{eq:optimprob2}
  {\mathcal C}_\Omega(\psi)
  =
  1 - \frac{\bar\bb^t \bD \bar\bb}{\|\bar\bb\|^2}
\end{equation}
and the minimization of
(\ref{eq:optimprob1}) becomes an eigenvalue problem.
The solution of this minimization depends on the opening
$\theta_0$. In Figure~\ref{fig:prol_local} we plot the value of
${\mathcal C}_{\Omega_{\theta_0}}$ against $\theta_0$ for
$\psi^\star_{\Omega_{1^{\circ}}},\psi^\star_{\Omega_{5^{\circ}}},\psi^\star_{\Omega_{10^{\circ}}}$.
The lowest curve is the minimum of the criterion for all openings
$\theta_0$. It is clear that there is no optimal function uniformly
in $\theta_0$: the concentration criterion $\mathcal{C}_{\Omega_0}$ of
each PSWF $\psi^\star_{\Omega_{\theta_1}}$ reaches the best
possible value for $\theta_0=\theta_1$ only.

\begin{figure}[htbp]
  \centering\includegraphics{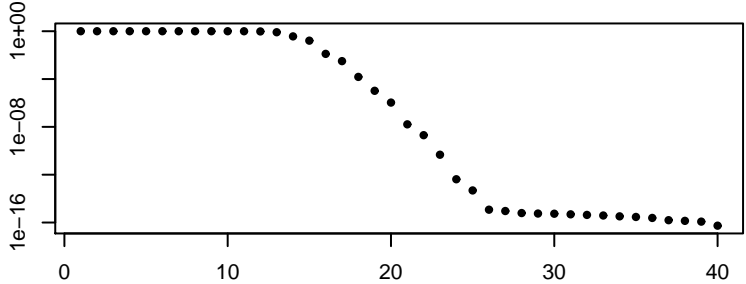}
  \caption{ \label{fig:D_eigenvalues} Eigenvalues of matrix D with
    $\theta_0 = 50^\circ$ and $L=[17,64]$. In this case, Shannon
    number $N = 13.3$.}
\end{figure}

As in the 1-dimensional case, the spectrum of $\bD$ exhibits a ``step
function'' behaviour: denoting $ N =\tr\bD $ (the ``Shannon
number''),
the matrix $\bD$ has about $N$ eigenvalues very close to 1, and most
of the others close to zero (see Fig.\ref{fig:D_eigenvalues},
and \citealt{Simons+2006} for details).

When several eigenvalues of $\bD$ are extremely close to 1, it is
computationally difficult to find the largest one and the associated
eigenvector.  In the case of $\Omega$ a polar cap and $\lmin = 0$, one
can advantageously solve the less degenerated eigenvalue problem
associated with the Gr\"unbaum differential
equation~\citep{grunbaum:longhi:perlstadt:1982} which has the same
solutions as~(\ref{eq:optimprob1}). We are not aware of an
equivalent theory in the case $\lmin > 0$.

\begin{figure}
  \centering
  \includegraphics{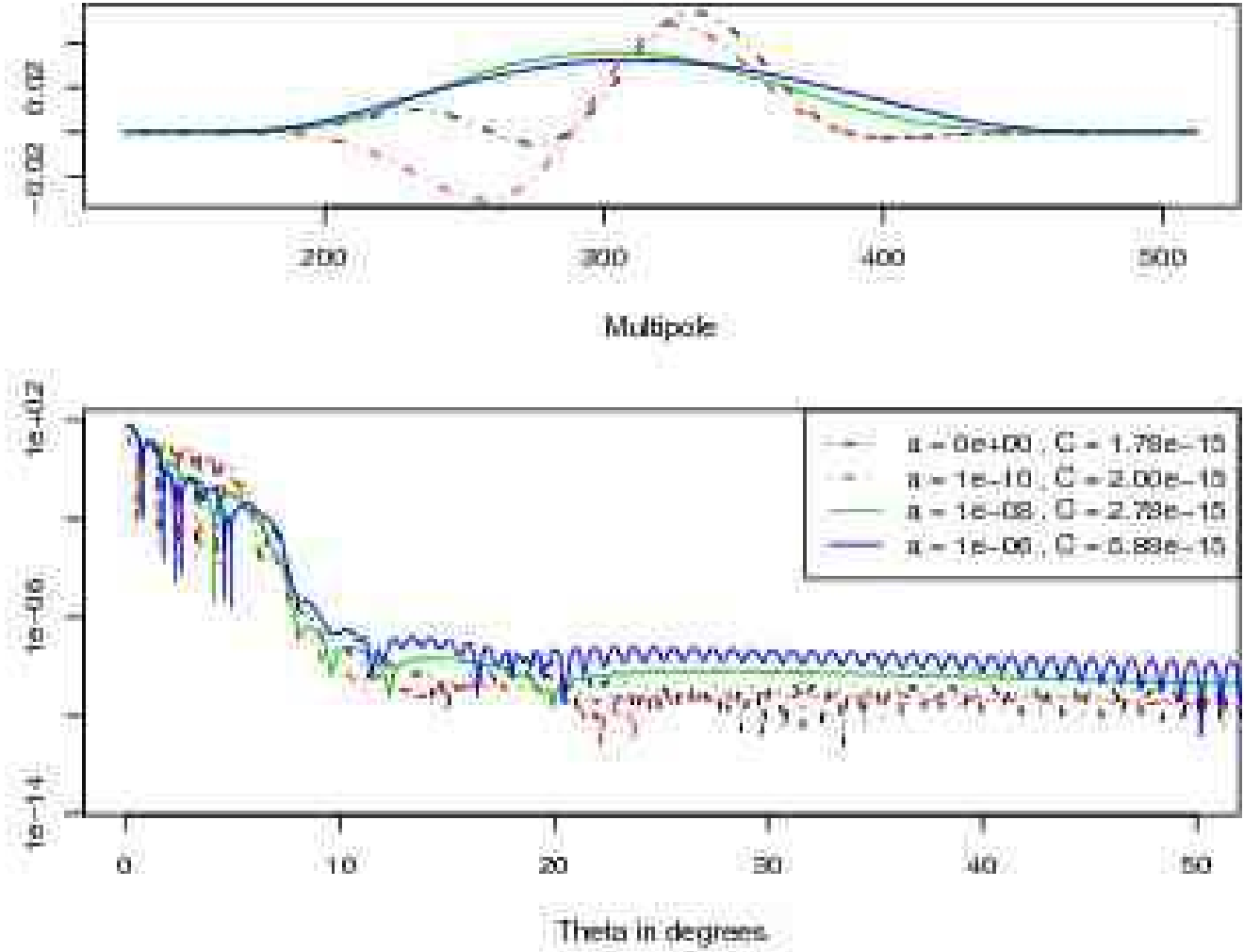}
  \caption{Effect of the smoothing on the spectral and spatial shapes
    of PSWFs.} \label{fig:smoothing}
\end{figure}

With $\epsilon$ being of the order the machine precision, all vectors
in $V_\epsilon=\bigoplus\limits_{\lambda\geq1-\epsilon}
\textrm{Ker}(\bD-\lambda\mathbf{Id}) $ have well spatially localized
counterparts, but they are not necessarily positive (in harmonic
domain).  This is not acceptable for instance if we were to use them as windows
associated to smoothing operator (denoted $h$ in the first Section),
and implement this operator using a needlet analysis-synthesis scheme, the
window of which has to be
defined as the square-root of the PSWF's window. 
To circumvent this,  we therefore introduce a modified coupling matrix $\widetilde{\bD} =
\bD+a\bH^t\bH $ where $a>0$ is a tuning parameter and $\bH$ is the
tridiagonal second-order finite difference matrix.
Window functions are now obtained as minimizers of
$\widetilde{\mathcal{C}_\Omega} (\psi) = 1 - \frac{\bar\bb^t
  \widetilde \bD \bar\bb}{\|\bar\bb\|^2}$ instead of
$\mathcal{C}_\Omega$.
The additional term favors non-oscillating functions among the vectors
of $V_\epsilon$ which are undistinguishable from their eigenvalues
$\lambda$.  Adding the ``smoothing'' term is  expected not to alter the
spatial localization of the filter. In practice, parameter $a$ is
selected to ensure `computational uniqueness' of the smallest
eigenvalue of $\widetilde \bD $.
Solutions obtained by the numerical implementation
of the minimization of $\widetilde{\mathcal{C}_\Omega}$ are displayed
in Figure~\ref{fig:smoothing}, with various values for the smoothing
parameter $a$. 
Dashed lines correspond to the vector returned numerically as the
``best'' eigenvector of $\bD$ (associated to the greatest eigenvalue),
and the best eigenvector of $\tilde \bD$ with parameter $a$ chosen
deliberately too small to ensure computationally uniqueness.
Oscillating functions are indeed obtained.
As $a$ grows, the criterion selects non oscillating windows, two of
which are shown by the plain lines.  The loss measured by the
increase of $\mathcal{C}_\Omega$ is displayed in the legend of the
lower panel and appears extremely small. In our example, the energy
outside $\Omega$ for the needlet built from $\tilde{\mathcal{C}}$
takes the value $2.78.10^{-15}$, whereas its minimal possible value is
$1.78.10^{-15}$.

A generalization of the Slepian concentration problem can be to
consider other measures of concentration, such as $\Lset^p$,
$p=1,...,\infty$ instead of $\Lset^2$. The criterion defined in
Eq.~(\ref{eq:optimprob1}) becomes then
\begin{equation}
  \label{eq:Cp}
  \mathcal{C}^{(p)}_\Omega(\psi)
  =
  1-\frac{\|\psi\ind_\Omega\|_p^p}{\|\psi\|_p^p}
\end{equation}
where $\|f\|_p^p=\int_\Sset|f(\xi)|^p\dd\xi$ if $p\in[1,\infty)$ and
$\|f\|_\infty^\infty=\textrm{ess}\sup\limits_{\xi\in\Sset}|f(\xi)|$
for a spherical function $f$.  Unlike Slepian criterion
$\mathcal{C}_\Omega=\mathcal{C}_\Omega^{(2)}$, these alternate
criteria do not lead to simple eigenvalue problems.  They could be
numerically optimized but this is beyond the scope of this paper.
However we compare in Section \ref{sec:table_comp} this criterion to
the original one $\mathcal C_\Omega$.

\subsection{Statistical criterion for optimal analysis with missing data}
\label{sec:statistical-criteria}

Instead of focusing on the ``geometrical'' shape of the needlet, one
may also optimize directly some alternate criterion of practical
interest.  

In this section, we consider the following framework: given an
underlying random field $X$ on $\Sset$ to be analysed, a window
function $W$ on $\Sset$ multiplying the field (for example a mask
putting the field to zero in some regions) and a region
$\Dset\subset\Sset$ of interest in which the analysis is to be done,
the aim is to get, in $\Dset$, needlet coefficients of $WX$ as close
as possible to the coefficients computed from the uncorrupted field
$X$.

We shall assume statistical properties on the fields $X, W, D$ and look for optimality of the
filters on average.
\begin{hyp}
\label{hyp:gaussian}
\begin{enumerate}
\item $X$ is a real-valued Gaussian zero mean isotropic square integrable
  random field on $\Sset$, with power spectrum $(C_\ell)$.
\item   $W$ and $D$ are deterministic elements of $\Hset$. 
\end{enumerate}
\end{hyp}
Implicitly, $X$ is a measurable mapping from some $(\mathcal X,
\mathsf X, \proba )$ into $(\Hset,\mathsf H)$, $\mathsf H$ being the
Borel $\sigma$-filed of $\Hset$.  Let $\esp$ denote the expectation
operator under $\proba$.  Recall that under
Assumption~\ref{hyp:gaussian}, the covariance function on the field
$X$ is well defined and is given by
\begin{equation*}
  \esp[X(\xi) X(\xi')] = (4\pi)^{-1}\sum_{\ell \in \Nset}  C_\ell L_\ell(\xi
  \cdot \xi') \ .
\end{equation*}
It follows that $ \esp[X(\xi)^2] = (4\pi)^{-1}\sum_{\ell \in \Nset}
(2\ell + 1) C_\ell $.  Moreover, the multipole moments $(a\lm)$ of $X$
are complex Gaussian random variables. They are centered, independent up to the relation
$a\lm=a^*_{\ell,-m}$ and satisfy $\esp(|a_{\ell
  0}|^2)=\frac12\esp(|a\lm|^2)=C_\ell$, $m\neq0$.

Note that $W$ and $D$ can be indicator functions (binary \emph{masks}) or any
smooth functions on the sphere.

A first attempt in this direction is the derivation of an unbiased
estimate of the spectrum from the multipole moments and the empirical
power spectrum of the weighted sky $XW$ defined by $\hat
a\lm=\int_\Sset X(\xi) W(\xi) Y\lm^*(\xi) \dd\xi$ and
$\hat{C}_\ell=\frac{1}{2\ell+1}\sum\limits_m\hat{a}\lm^2$
respectively. It is well-known \citep[see][see also the compact proof
in Appendix \ref{sec:proofs}]{Peebles:1973, Hivon+2002}
that
\begin{equation}
  \label{eq:master}
\esp(\hat{C}_\ell)=\sum_{\ell' \in \Nset}\mathcal M_{\ell
  \ell'}C_{\ell'}\text{ with }\mathcal M_{\ell \ell'} =
\sum\limits_{0\leq\ell''\leq\ell+\ell'}\alpha_{\ell \ell'
  \ell''}\frac{2\ell''+1}{2\ell+1}C_{\ell''}^W \; ,
\end{equation}
where the coefficients $\alpha_{\ell \ell' \ell''}$ are defined by
(\ref{eq:defalpha}).  Note that the coupling matrix $\mathcal M$
depends on $W$ only through its `power spectrum' $C_\ell^W$.  If
$\mathcal{M}$ is invertible, then $(\mathcal M^{-1}(\hat
C_{\ell'}))$ provides an unbiased estimate of $(C_\ell)$.

Let now derive a criterion to design a window function $b$ which
minimises the effect of missing data in a needlet analysis procedure.
We focus on a single band smoothed field (\ie{} \ we fix one scale
$j$) and the dependence on $j$ is implicit in the notations. For a
collection of couple of indices, say $(\ell_i,m_i)_{i=1,\dots,I}$, we
use $\sumstar_{(\ell_i,m_i)_{i=1,\cdots,I}}$ as a shorthand notation
for the summation on $\ell_i \in \Nset, m_i \in
\{-\ell_i,\cdots,\ell_i\}, i=1,\dots,I$.

Given an analysis spectral window $\bb =
(b_{\ell_{\min}},\cdots,b_{\ell_{\max}})$ and its associated smoothing
operator $\Phi=\sum_{\lmin\leq\ell\leq\lmax}b_\ell\Pi_\ell$, the
smoothed masked field is $$\Phi XW(\xi) = \sum_{\ell \in L} b_\ell
\int_\Sset X(\xi')W(\xi')L_\ell(\xi\cdot\xi')\dd\xi'.$$
Write $\esp[\Phi X(\xi)^2] = (4\pi)^{-1}\sum_{\ell} \sigma_\ell^2
b_\ell^2$ with $\sigma_\ell^2 = (2\ell + 1) C_\ell$.
Let $\epsilon$
denote the normalized difference field
\begin{align}
  \epsilon(\xi) &= \frac{\Phi X(\xi) -
    \Phi (XW)(\xi)}{\esp^{1/2}[\Phi X(\xi)^2]} \nonumber  \\
  &= \left(\sum \sigma_\ell^2b_\ell^2\right)^{-1/2} \sumstar_{(l,m)} b_\ell
  \bar a\lm Y\lm(\xi)
\label{eq:expepsilon}
\end{align}
where we have defined $\bar W = 1 - W$, $\bar a\lm = \langle X\bar W,Y\lm
\rangle$.

Suppose that $(\{(\xi_k\}_{k \in K},\{\lambda_k)\}_{k \in K})$ provides
an exact Gauss quadrature formula at a degree $2\lmax$.  Define
$\beta_k$ and $\beta'_k$ the needlet coefficients of $X$ and
$XW$, respectively and define $$\epsilon_k = \frac{\beta_k -
  \beta'_k}{\sqrt{\esp(\beta_k^2)}}.$$  Those random variables are
normalized errors on the needlet coefficients induced by the
application of the weight function $W$. If both $X$ and $XW$ are in
$\Hset_{\lmax}$, we easily check that $\esp(\beta_k^2) =
\sqrt{\lambda_k} (4\pi)^{-1} \sum_\ell b_\ell^2 \sigma^2_\ell$ and
\begin{equation*}
  \forall k \in K, \ \epsilon_k = \epsilon(\xi_k) \ .
\end{equation*}

The dispersion of either the continuous field $\epsilon(\xi)$ or the finite
set $\{\epsilon_k\}_{k \in K}$ is taken as a measure of quality for an
analysis $\Phi$.
This dispersion is not measured on the whole sphere, since the
difference $\epsilon$ must be important in the regions where $W$ is
far from~1.  In order to select the regions where $\epsilon$ is to be
minimized we introduce a function $D = \sum d\lm Y\lm$ which provides
a positive weight function in $\Hset$.
In the simplest case $D$ can be $\ind_\Dset$ for a region $\Dset$ of
interest. More generally, $D$ can be designed to give more or less
importance to various regions of $\Sset$ according, for instance, to
the need for reliability in the needlet coefficients.

The coefficients $\epsilon_k$ or their continuous version $\epsilon$
are used in two ways.  
The first one introduces a ``tolerance'' threshold $\alpha$ and counts
the number of coefficients which are on average below this threshold.
This measure of the efficiency of a filter in the presence of a mask
is presented in \citealt{baldi:etal:2006a,
  pietrobon:balbi:marinucci:2006} but its optimization was not
considered.
The second one considers the integrated square error of $\epsilon$,
weighted by the function $D$.  It leads to a quadratic quadratic which
is readily optimized.

The first criterion, writes, for a binary function $D$,
\begin{equation}\label{eq:mask_error}
  E_\bb(\alpha)= \frac{\sum_{k : D(\xi_k)=1} \proba(|\epsilon_k| <
    \alpha)}{\sharp\{k:D(\xi_k)=1\}},
\end{equation}
that is, the mean fraction of needlet coefficients corrupted by less
than a normalized error $\alpha\geq0$.  For an arbitrary function $D$,
a possible generalization of~(\ref{eq:mask_error}) is
\begin{equation*}
  E_\bb(\alpha) = \frac{\sum_{k \in K} D(\xi_k)
    \proba(|\epsilon_k|\leq\alpha)}{\sum_{k \in K} D(\xi_k)}. 
\end{equation*}
In Subsection~\ref{sec:robustn-with-resp}, we compare different
windows using this criterion and a real mask.

Alternately, consider now the mean integrated square error (MISE)
 \begin{equation}
  \label{eq:defMISE1}
  R(\bb) = \esp \int_\Sset  D(\xi) \|\epsilon(\xi)\|^2 \dd \xi
\end{equation}
and define the optimal shape for the window $\bb$ as
\begin{equation}
\label{eq:defopt}  \bb^\star = \arg\min_{\|\bb\| = 1 }  R(\bb).
\end{equation}
Straightforward algebra leads to a close form expression of $R(\bb)$
depending on $\bb$, on the weight functions $W$ and $D$, and on the
power spectrum $(C_\ell)_{\ell \in \Nset}$.  Let $\bar w\lm, d\lm$
denote the multipole coefficients of the weight functions $\bar W, D$,
respectively and
$$\varwigner{\ell}{\ell'}{\ell''}{m}{m'}{m''} :=  \int_\Sset Y\lm
(\xi)Y_{\ell'm'}(\xi)Y_{\ell''m''}^*(\xi)\dd\xi$$ (see~(\ref{eq:wigner3j1}) for an
expression as a function of the Wigner-3$j$ coefficients).
\begin{prop}\label{prop:expstatcriterion}
Under Assumption~\ref{hyp:gaussian}
\begin{equation*}
  R(\bb) = \frac{\bb' \bQ
  \bb}{\bb' \boldsymbol \sigma \bb}
\end{equation*}
where $\boldsymbol \sigma = \mathrm{diag}((\sigma_\ell^2))$ and $\bQ$ is the
matrix with entries
\begin{multline}
  Q_{\ell \ell'} = \sum_{m,m'}\sumstar_{(\ell_1,m_1)}C_{\ell_1} \\ \sumstar_{(\ell_i,m_i)_{i=2,3,4}}
   \bar w_{\ell_2 m_2} \bar w_{\ell_3 m_3}^* d_{\ell_4 m_4} 
 \varwigner{\ell_1}{\ell_2}{\ell}{m_1}{m_2}{m}
 \varwigner{\ell_1}{\ell_3}{\ell'}{m_1}{m_3}{m'}^*
 \varwigner{\ell}{\ell_4}{\ell'}{m}{m_4}{m'} \label{eq:Qllprime} \ .
\end{multline}  
If both $W$ and $D$ are axisymmetric the ten-tuple summations above
reduce to a five-tuple one
\begin{align*}
  Q_{\ell \ell'} &= \sum_{m}  
\sum_{\ell_1,\ell_2,\ell_3,\ell_4} C_{\ell_1} \bar w_{\ell_2 0}
  \bar w_{\ell_3 0} d_{\ell_4,0} \varwigner{\ell_1}{\ell_2}{\ell}{m}{0}{m}
  \varwigner{\ell_1}{\ell_3}{\ell'}{m}{0}{m}
  \varwigner{\ell}{\ell_4}{\ell'}{m}{0}{m}\\ &= \sum_{m} A_{\ell \ell' m } D_{\ell \ell' m}
\label{eq:Qllprimesimple} \ .
\end{align*}
\end{prop}

In the next section we shall give some illustrative examples of optimal
spectral windows $\bh^\star$ in the
particular axisymmetric case.
\begin{rem}
  As in the Slepian's problem, the design of an optimal filter reduces
  to an eigenvalue problem.  In particular, if $\sigma_\ell >0$ for
  any $\ell \in L$, write $ b^\dagger_\ell = \sigma_\ell b_\ell$. Let
  $\bb^{\dagger \star}$ be an eigenvector associated with the lowest
  eigenvalue of $\bQ^{\dagger}$, $Q_{\ell \ell'}^{\dagger} =
  (\sigma_\ell\sigma_{\ell'})^{-1} Q_{\ell \ell'}$. Then $\bb^\star :=
  \boldsymbol \sigma \tilde \bb^{\dagger \star} / \| \boldsymbol \sigma \tilde
  \bb^{\dagger \star} \|$ is a solution of~(\ref{eq:defopt}).
\end{rem}
\begin{rem}
  \label{rem:sums}
  For those sums to be tractable, one has to assume that $D$, $W$,
  $C_\ell$ have finite support in the frequency domain, \ie{} that the
  windows $D$ and $W$ are smooth (or apodized) and $C_\ell = 0$ for
  large enough $\ell$.
\end{rem}
\begin{rem}
  \label{rem:truemask}
  The matrix $\bQ$ being a second-order moment for the random field
  $X$, it can also be approximated by a moment estimator using
  Monte-Carlo experiments. This remark is of important practical
  interest as we are mostly concerned with non zonal masks.
\end{rem}

\section{Examples, numerical results}
\label{sec:results}

\subsection{Comparison of filters for various criteria}
\label{sec:table_comp}

In Section~\ref{sec:crit-filt-design}, we considered several criteria
measuring the localization properties of filters, and derived explicit
or computational optimization for some of them. In
Table~\ref{tab:variouscrit}, we compare the scores reached by the
filters displayed in Figure~\ref{fig:examples}.
The columns indexed by $\Lset^2$-$\theta$ list the values
$\mathcal{C}_{\Omega_\theta}(\psi)$ defined in
Eq.~(\ref{eq:optimprob1}).  More generally, the columns indexed by
$\Lset^p$-$\theta$ correspond to the values
$\mathcal{C}^p_{\Omega_\theta}(\psi)$ defined in Eq.~(\ref{eq:Cp}). A
column lists the values of $1-E(\alpha)$ defined in
Eq.~(\ref{eq:mask_error}), applied with the mask Kp0 of
Fig.~\ref{fig:WMask} and a tolerance parameter $\alpha=10\%$ (see next
subsection for more details). A last column gives, by way of illustration
only, the value of the ``uncertainty product''
$\Delta_\xi(\psi)\times\Delta_{\mathsf{L}}(\psi)$, where
\begin{equation}
  \label{eq:Heisenberg}
  \Delta_\xi(\psi) =
  \frac{\sqrt{1-\|\int_\Sset\xi\psi(\xi)^2\dd\xi\|^2}}{\int_\Sset\xi\psi(\xi)^2\dd\xi}
\text{ and } \Delta_{\mathsf{L}}(\psi) = \sum_{\ell\geq0}\ell(\ell+1)b_\ell^2.
\end{equation}
\cite{Narcowich+96} proved that
$\Delta_\xi(\psi)\times\Delta_{\mathsf{L}}(\psi)\geq1.$

\begin{sidewaystable}
  \centering
  \scriptsize
  \begin{tabular}[t]{|c|c|c|c|c|c|c|c|c|c|c|c|c|c|c|}
      \hline  
& $\Lset^2$-0.5$^\circ$ & $\Lset^2$-1$^\circ$ & $\Lset^2$-1.5$^\circ$ & $\Lset^2$-5$^\circ$ 
& $\Lset^1$-0.5$^\circ$ & $\Lset^1$-1$^\circ$ & $\Lset^1$-1.5$^\circ$ & $\Lset^1$-5$^\circ$ 
& $\Lset^\infty$-0.5$^\circ$ & $\Lset^\infty$-1$^\circ$ & $\Lset^\infty$-1.5$^\circ$ & $\Lset^\infty$-5$^\circ$ 
& 1-$E$(0.1) & $\Delta_\xi\Delta_\mathsf{L}$\\
\hline\hline

Spline, order 3 &2.2e-02 & 5.2e-03 &  7.4e-04 &  9.8e-07 &  4.2e-01 &
2.2e-01 & 1.0e-01 & 1.5e-02 & 5.0e-02 & 1.9e-02 & 5.1e-03 & 6.4e-05 &
2.6e-01 & 2.7 \\
\hline
Spline, order 7 &4.0e-02 & 1.3e-02 &  2.0e-03 &  4.8e-08 &  5.0e-01 &
2.9e-01 & 1.3e-01 & 1.7e-03 & 6.0e-02 & 2.7e-02 & 7.1e-03 & 1.2e-05 &
3.3e-01 & 3.1 \\
\hline
Spline, order 15 &6.1e-02 & 2.5e-02 &  4.9e-03 &  4.0e-07 &  5.9e-01 &
4.0e-01 & 2.2e-01 & 2.3e-03 & 6.9e-02 & 3.3e-02 & 9.8e-03 & 7.0e-05 &
4.1e-01 & 3.7 \\
\hline
Spline, order 21 &7.2e-02 & 3.1e-02 &  7.1e-03 &  7.7e-06 &  6.2e-01 &
4.5e-01 & 2.7e-01 & 1.0e-02 & 7.3e-02 & 3.7e-02 & 1.1e-02 & 2.7e-04 &
4.6e-01 & 4.1 \\
\hline
Prolate, cap 0.5$^\circ$ &\textcolor{blue}{\bf1.2e-02} & 6.0e-03 &
3.4e-03 &  9.5e-04 &  8.5e-01 &  8.2e-01 & 8.0e-01 & 7.2e-01 & 5.1e-02
& 1.0e-02 & 5.6e-03 & 1.0e-03 & 6.5e-01 & 9.8 \\
\hline
Prolate, cap 1$^\circ$ & 6.7e-02 & \textcolor{blue}{\bf4.3e-05} &
5.8e-06 &  1.7e-06 &  \textcolor{blue}{\bf3.8e-01} &  1.3e-01 & 1.2e-01 &
1.1e-01 & 1.1e-01 & \textcolor{blue}{\bf2.0e-03} & 2.0e-04 & 5.0e-05 &
\textcolor{blue}{\bf1.5e-01} & 3.1 \\
\hline
Prolate, cap 1.5$^\circ$ & 1.2e-01 & 1.5e-03 &
\textcolor{blue}{\bf3.4e-07} &  1.2e-08 &  4.3e-01 &
\textcolor{blue}{\bf5.3e-02} & \textcolor{blue}{\bf1.0e-02} & 8.8e-03 &
1.3e-01 & 1.7e-02 & \textcolor{blue}{\bf1.4e-04} & 4.5e-06 & 1.7e-01 &
3.6 \\
\hline
Prolate, cap 5$^\circ$  &1.1e-01 & 6.7e-03 &  6.5e-04 &
\textcolor{blue}{\bf5.7e-14} &  5.0e-01 &  1.8e-01 & 6.8e-02 &
\textcolor{blue}{\bf1.1e-06} & 1.2e-01 & 2.2e-02 & 5.9e-03 &
\textcolor{blue}{\bf2.6e-08} & 2.4e-01 & 3.6 \\
\hline
Exponential &1.8e-02 & 3.2e-03 &  1.0e-03 &  1.0e-05 &  4.4e-01 &
2.6e-01 & 1.9e-01 & 4.8e-02 & \textcolor{blue}{\bf4.4e-02} & 1.4e-02 &
5.7e-03 & 1.9e-04 & 2.7e-01 & \textcolor{blue}{\bf2.7} \\
\hline\hline
B-Spline &\textcolor{blue}{\bf1.1e-02} & 1.3e-03 &  3.9e-04 &  1.3e-05 &
4.8e-01 &  3.3e-01 & 2.7e-01 & 1.5e-01 & \textcolor{blue}{\bf3.1e-02} &
6.8e-03 & 2.5e-03 & 1.5e-04 & 2.1e-01 & \textcolor{blue}{\bf1.2} \\
\hline
Mexican hat &6.4e-01 & 1.1e-02 &  8.5e-07 &  7.3e-12 &  7.9e-01 &
8.8e-02 & \textcolor{blue}{\bf7.1e-04} & 1.6e-04 & 4.7e-01 & 7.9e-02 &
4.5e-04 & 1.4e-07 & 4.9e-01 & 3.0 \\
\hline

    \end{tabular}
    \caption{Comparison of the eleven filters of
      Fig.~\ref{fig:examples}, the nine first of which are
      band-limited in $L$=[256,1024].}\label{tab:variouscrit}
\end{sidewaystable}

The PSWFs perform the best not only for the $\Lset^2$ criterion which
they optimize, but also in most cases for the criteria where the
$\Lset^2$ norm is replaced by $\Lset^p$ ones, $p=1$ and $p=\infty$,
with the same opening angles $\theta_0$. Although the Kp0 mask has
many small cut areas all over the sphere, most of the 11 filters
presented here allow to retain more than $60\%$ of the outside-mask
coefficients $\beta_k$ if a $10\%$ error due to the presence of the
mask is accepted. The performance w.r.t. this criterion goes up to
$85\%$ for the PSWF optimally concentrated in a cap of
$1^\circ$. However, the choice of arbitrary value of $\alpha$ has a
major impact on the ranking of the filters. This point is investigated
in the next subsection.

\subsection{Robustness of needlets coefficients}
\label{sec:robustn-with-resp}

In this Subsection, we illustrate the performances of various window
functions using the criterion~(\ref{eq:mask_error}).  We have run
$N=30$ Monte-Carlo experiments to estimate the numerator of
$E_\bb(\alpha)$.  The random fields $X$ are drawn using the
$(C_\ell)$-spectrum of the best-fitting model for the CMB estimated by
the WMAP team \citep{Hinshaw+2006}. The mask $W$ was chosen as Kp0,
displayed in Figure~\ref{fig:WMask}, which masks the galactic plane
and many point sources. The band is $L=[256,1024]$.

\setlength{\widthsky}{8cm}
\begin{figure}[htbp]
  \centering    \includegraphics[width=\widthsky]{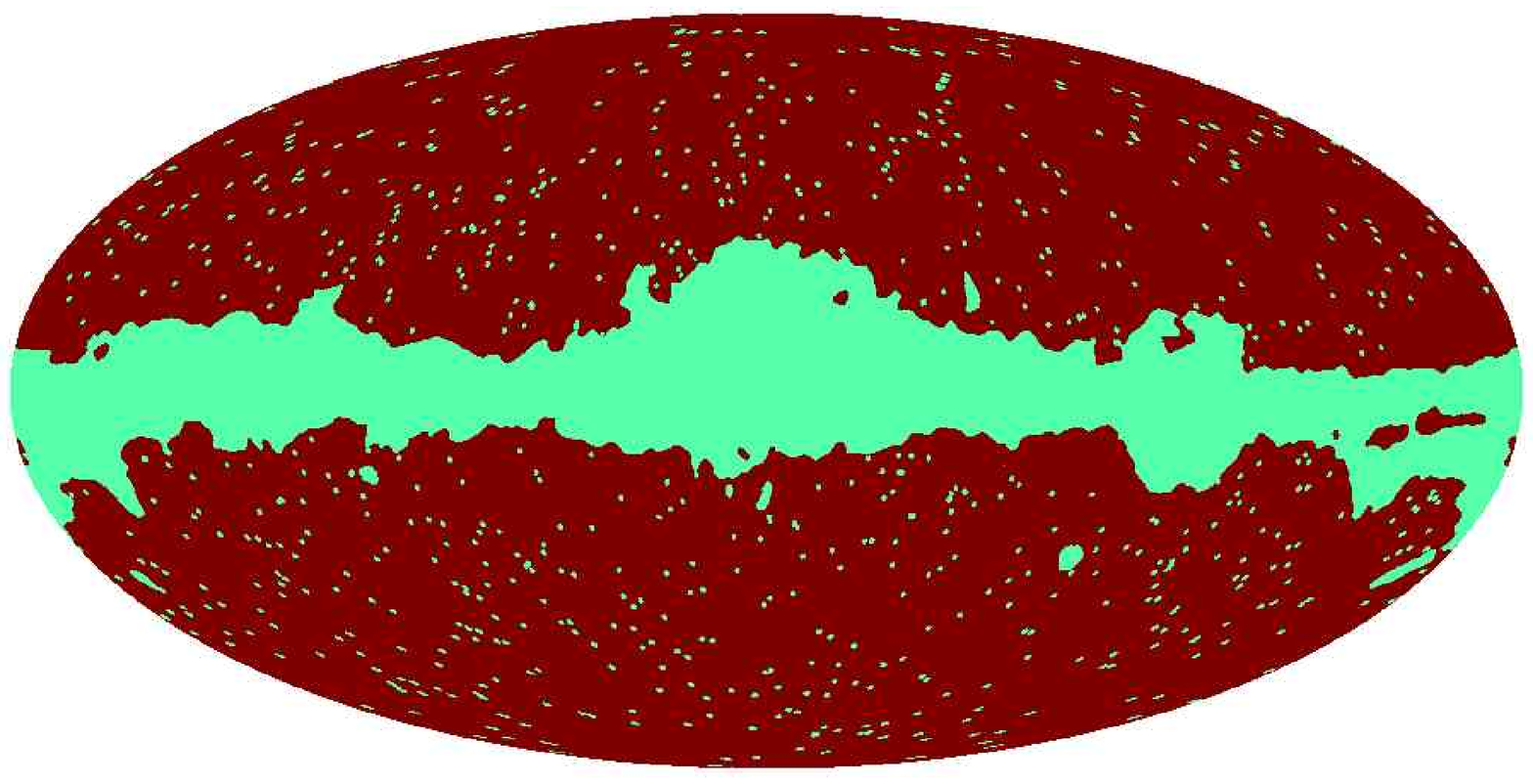}
  \caption{Kp0 mask.}  \label{fig:WMask}
\end{figure}

\begin{figure}[htbp]
  \centering    \includegraphics{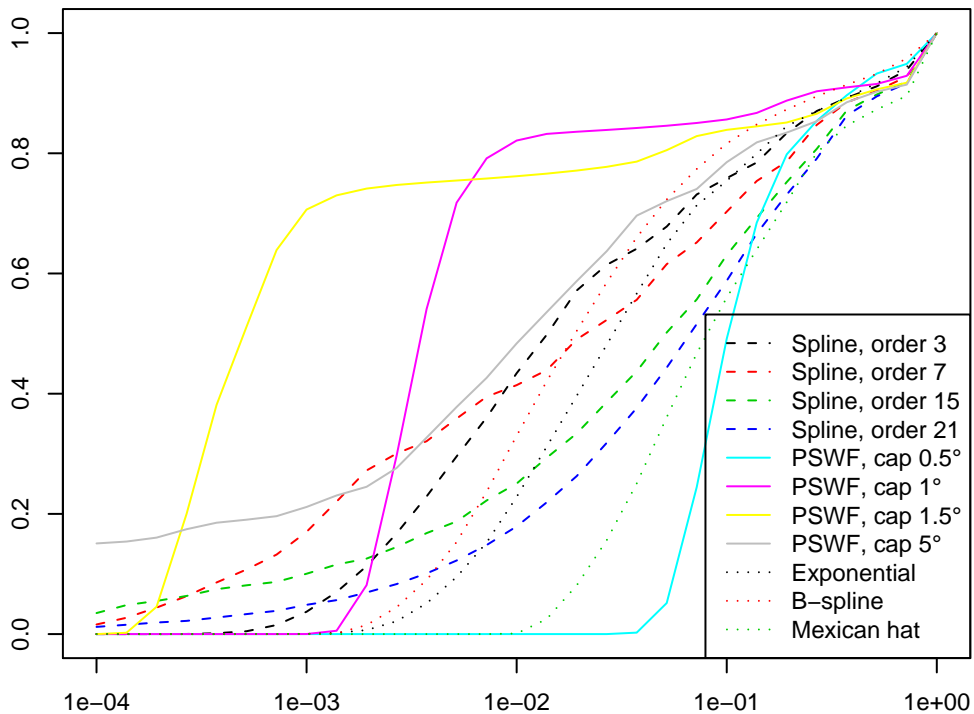}
  \caption{Proportion $E_\bb(\alpha)$ of coefficients uncontaminated
    at tolerance level $\alpha$.}  \label{fig:maskeffect}
\end{figure}

Figure~\ref{fig:maskeffect} compares the increasing functions
$E_\bb(\cdot)$ corresponding to various filters $\bb$.  There is no
``uniformly best'' (\ie{} highest in the figure) needlet: some allow
to retain more coefficients when the constraint imposed on the error
is loose enough, but their efficiency decreases faster as $\alpha$
goes to zero. Inspect \eg{} the PSWF family.

\setlength{\widthfig}{5cm}

\subsection{Some MISE-optimal filters for axisymmetric weight functions}

We present here the results of the optimization (\ref{eq:defopt}) in
the case of axisymmetric weight functions $W$. For simplicity, the
reconstruction weight function $D$ is taken equal to $W$.  We stick to
the CMB spectrum of previous subsection.

Figure~\ref{fig:masks} displays some of the masks $W$ used in the
experiments.  The apodization in simply a cosine-arch junction between
0 and 1, on a 2-degrees angular range. This means that the data is
available on the dark regions, and that its
$L=[\ell_{\min},\ell_{\max}]$-band-limited part has to be recovered in
this area too.  \setlength{\widthsky}{5cm}
\begin{figure}
  \centering
  \subfigure[]{\includegraphics[width=\widthsky]{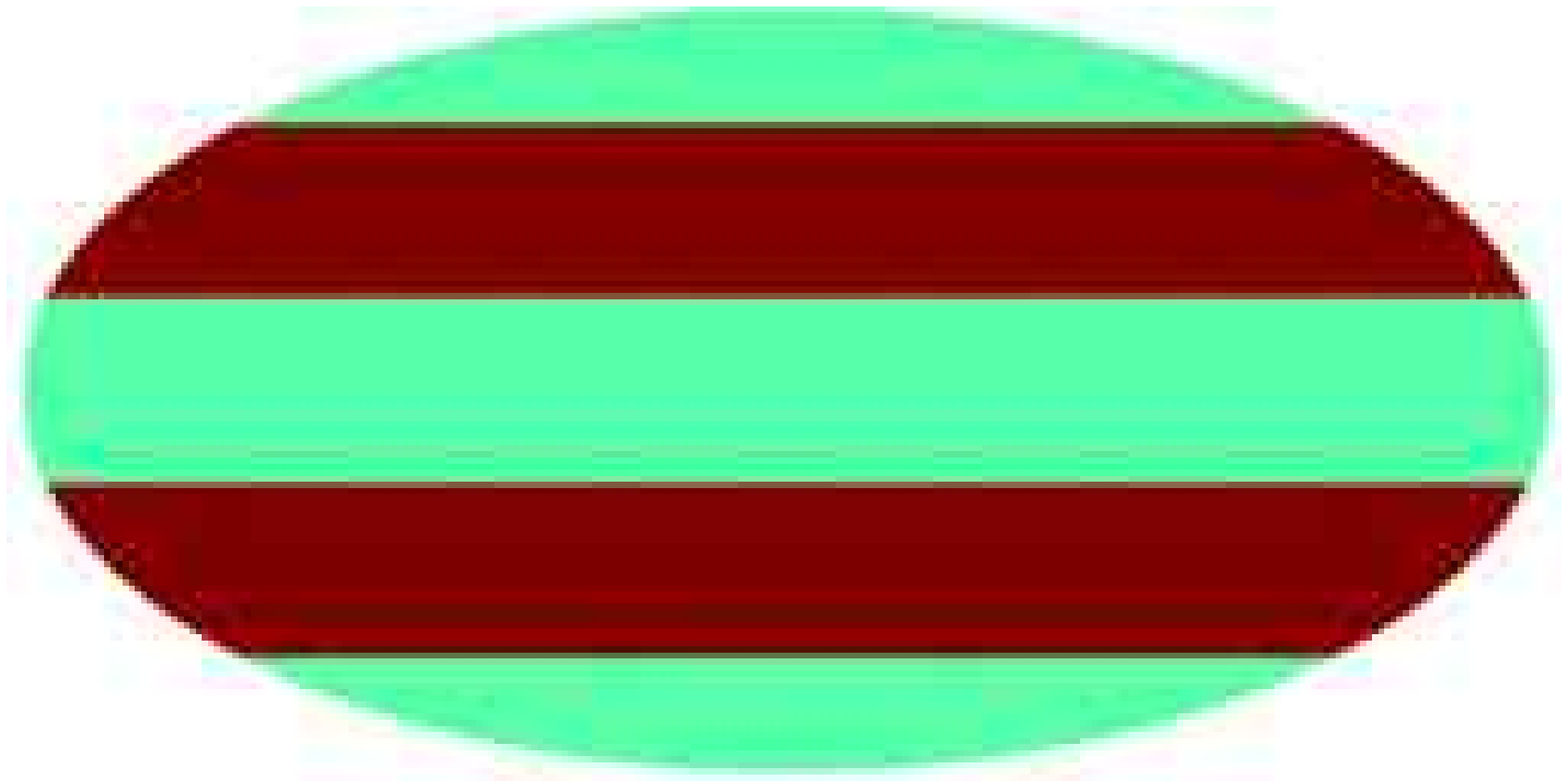}}
  \subfigure[]{\includegraphics[width=\widthsky]{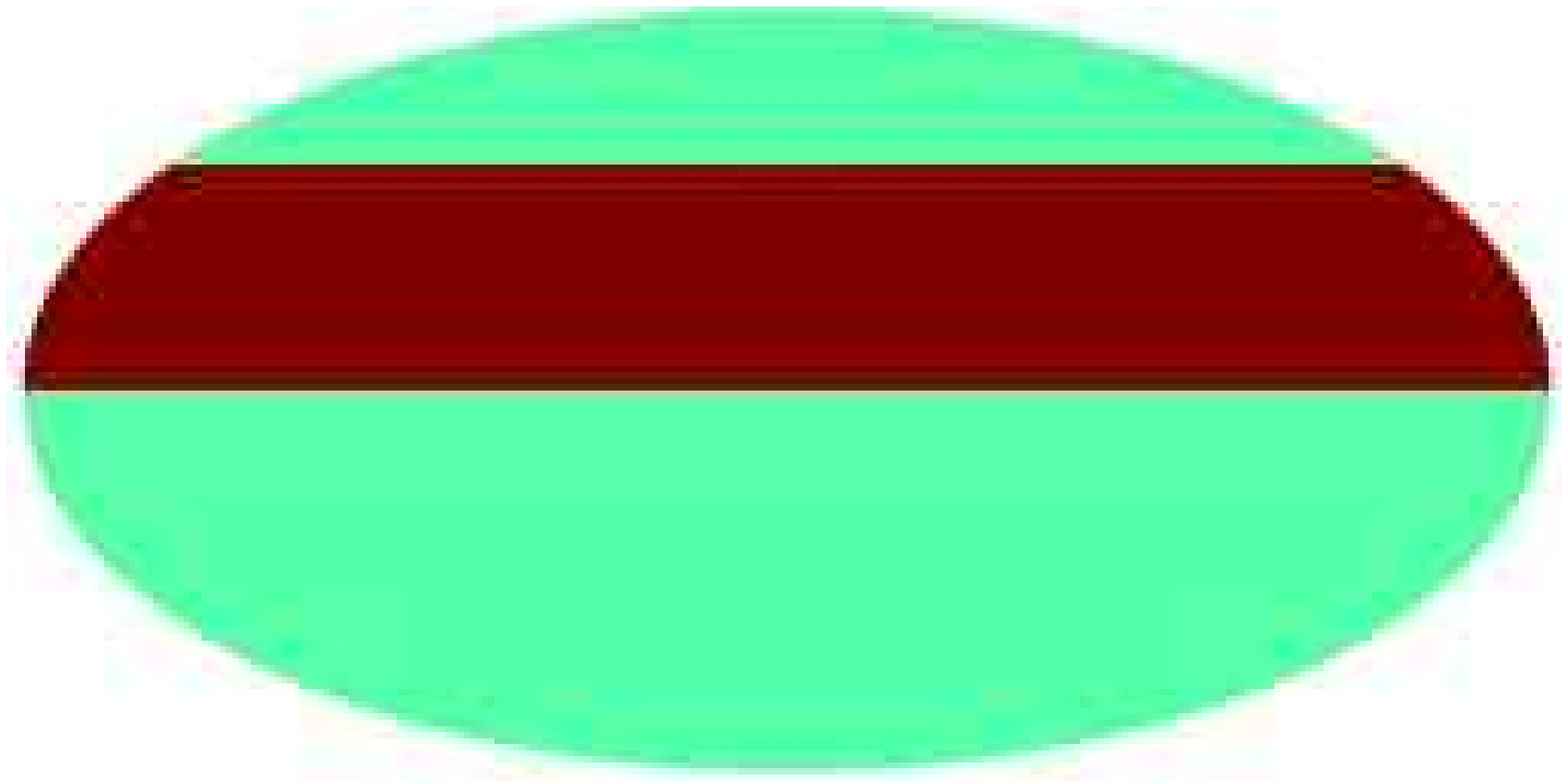}}
  \subfigure[]{\includegraphics[width=\widthsky]{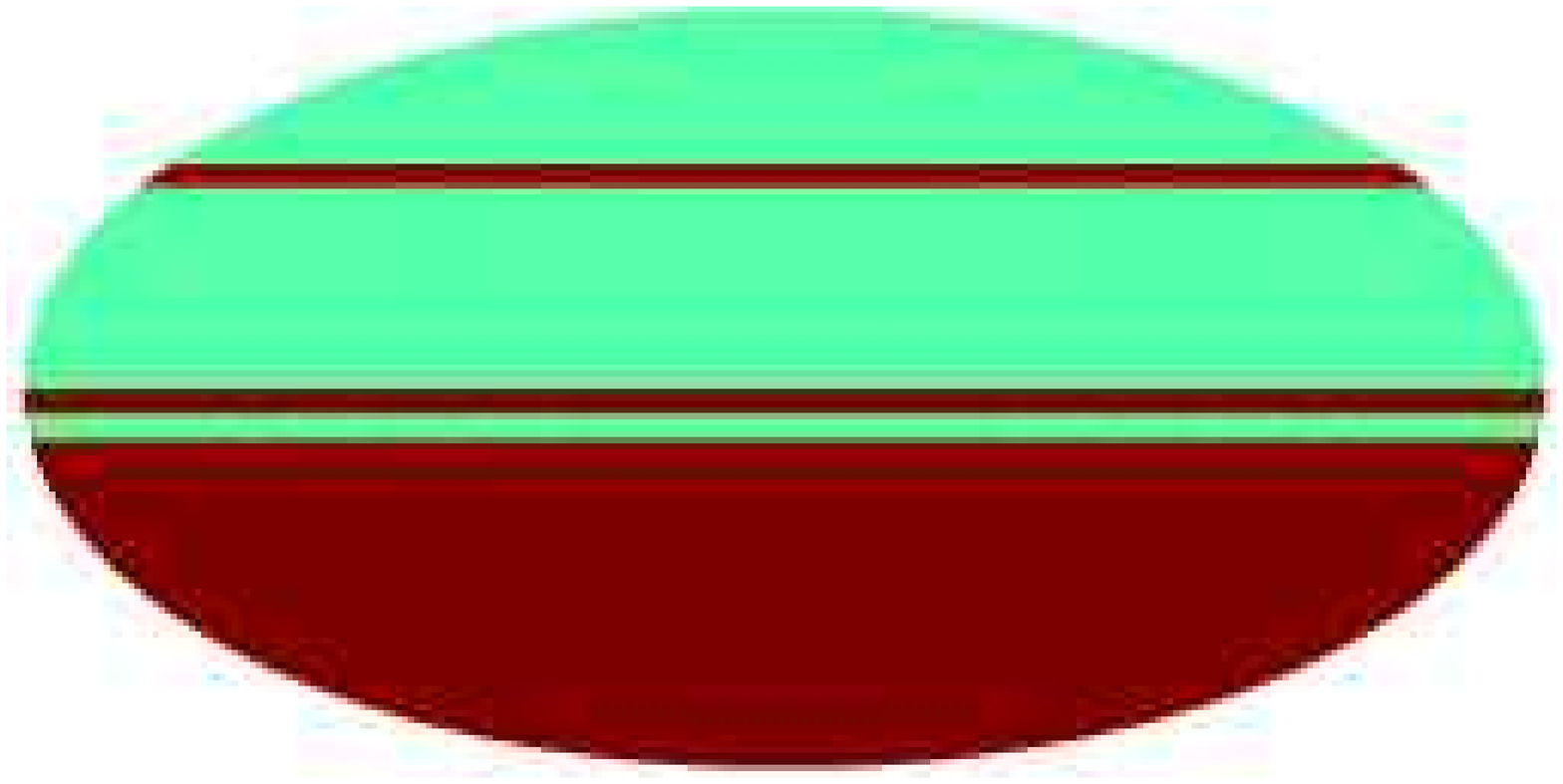}}
  \subfigure[]{\includegraphics[width=\widthsky]{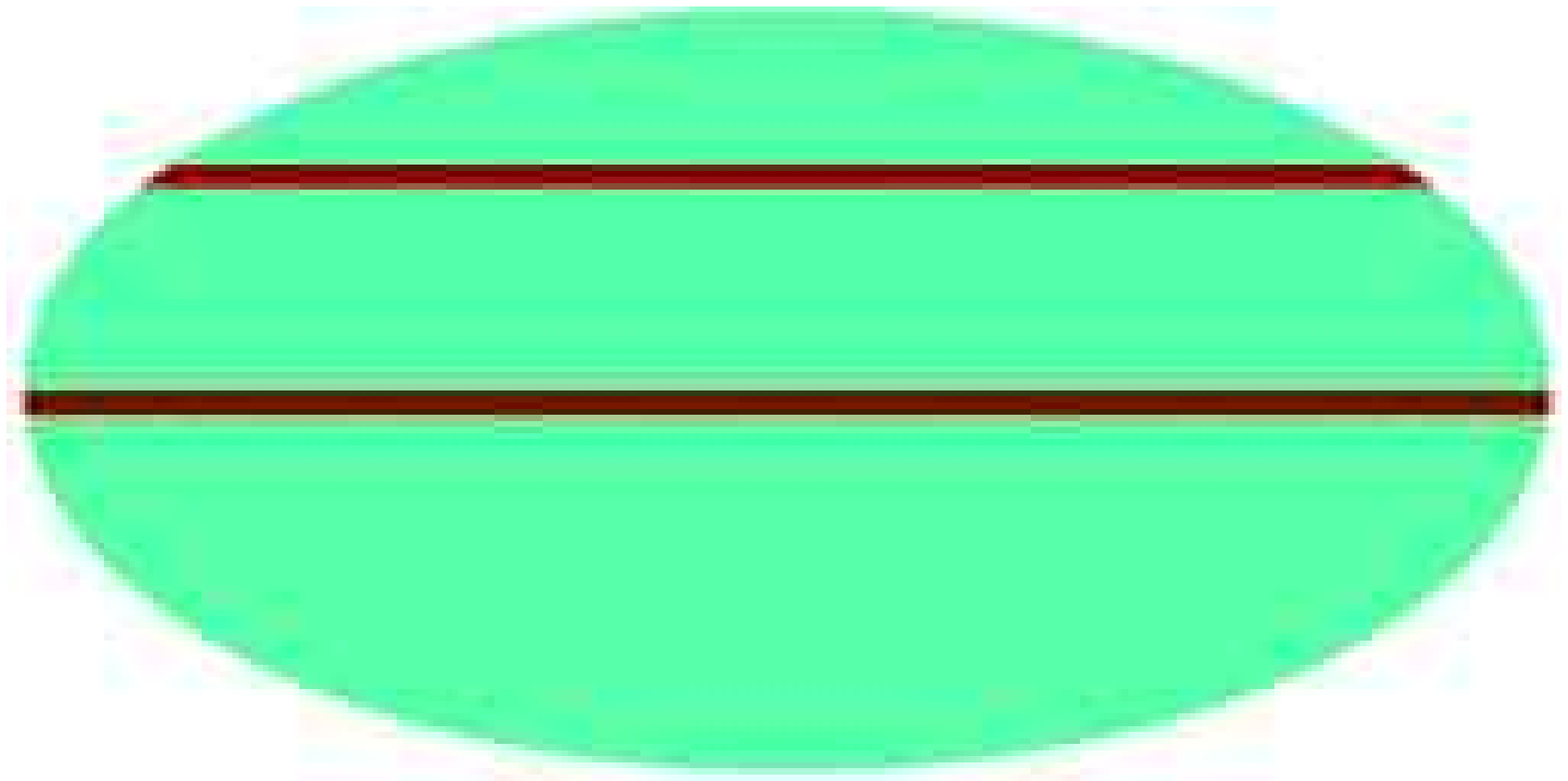}}
  \caption{Four different apodized masks. The degree of apodization, measured as the width of
    the cosine-arch 0-1 junction is, 2 degrees.}
  \label{fig:masks}
\end{figure}

\begin{figure}
  \centering
\includegraphics{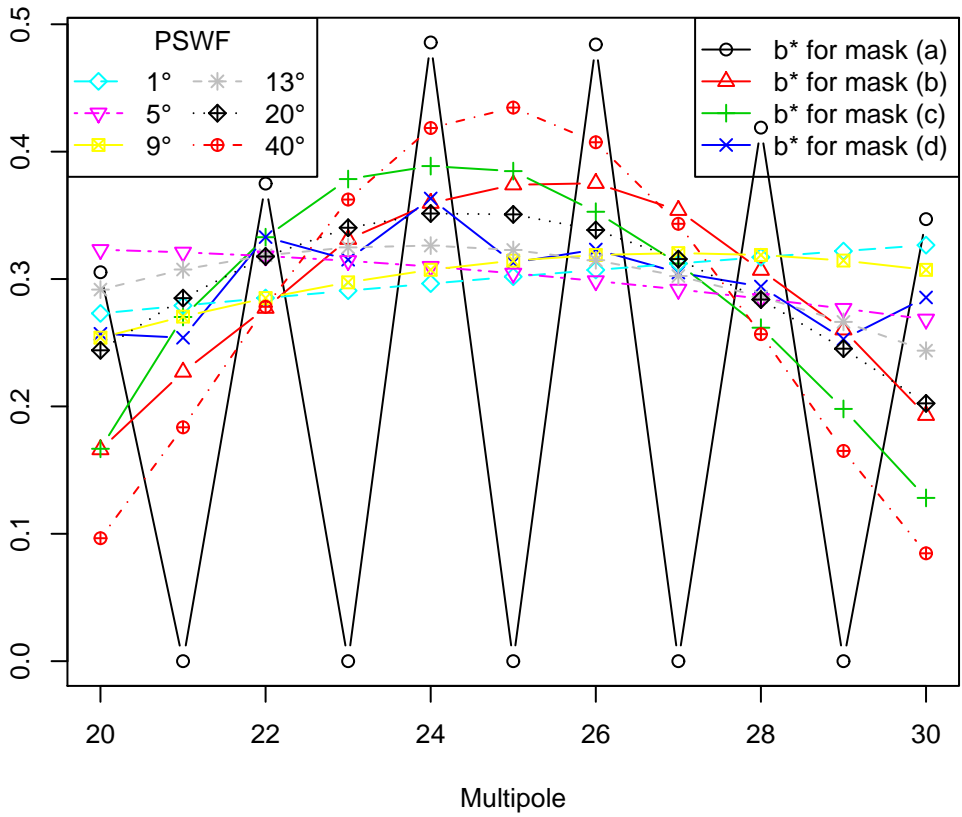}
    \caption{Shape of optimal window functions (plain lines) and PSWF
      (coloured and dashed lines) with various openings.}
  \label{fig:shape_optimal}
\end{figure}

On Figure~\ref{fig:shape_optimal} we have plotted the optimal filter in the
$R(\bb)$-sense for the masks of Figure~\ref{fig:masks} together with different
PSWFs.  The criterion captures the symmetry of the mask (a) (the shape of the
matrix $\bQ$ is a ``checkerboard''), and the optimal filter is thus zero on all
even (here) or all odd multipoles. The associated axisymmetric needlet $\psi$
is symmetric w.r.t. the equatorial plane, and thus is well concentrated around
both the North and the South poles. Such solutions are very sensitive to the
modifications of the masks.

\begin{figure}
  \centering
  \subfigure[]{\label{fig:bxpa}\includegraphics{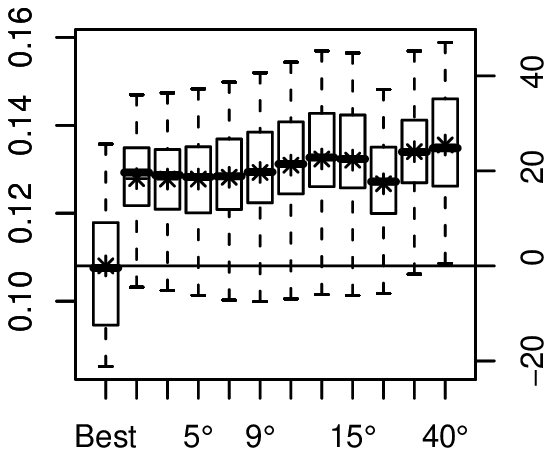}}
\subfigure[]{\label{fig:bxpb}\includegraphics{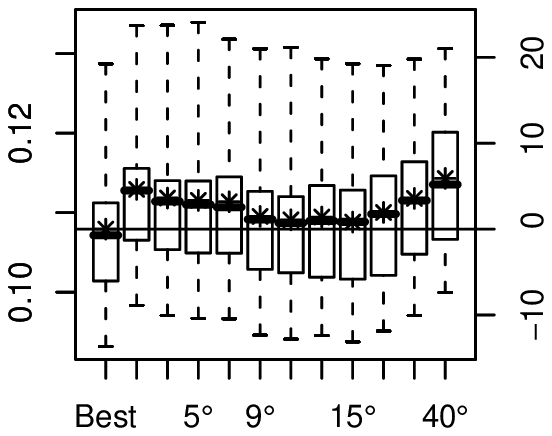}}
 \caption{Mean square error in analysis. Stars are potted at the
   estimated values for $R(\mathbf h)$. For Figure~\ref{fig:bxpa}, the mask is 
   Figure~\ref{fig:masks}(a) and $L = [5,15]$. For Figure~\ref{fig:bxpb}, the mask is
   Figure~\ref{fig:masks} and $L = [20,30]$}
  \label{fig:mse}
\end{figure}

We conducted a small Monte-Carlo study to confirm the benefit of our
approach. We have compared our best filters $\bb^\star$ to PSWFs
with different opening. On Figure~\ref{fig:mse}, we show the box-plots
of the distribution of the statistic $R(\mathbf h)$ for all those
filters. Stars are plotted at the position of the estimated value of
$\esp R(\bb)$ and the horizontal line is this value for
$\bb^\star$. The right vertical scale is for the relative error (in
percent) with respect to $\esp R(\bb^\star)$.

Fig.~\ref{fig:bxpa} illustrates the strong benefit of a filter that captures
the geometry of the mask. The relative improvement with respect to the best
PSWF is of order 20\%.  It should be noted however that the shape of this
optimal filter (described above) may lead to a misleading space-frequency
picture.  In some other cases, as shown in Figure~\ref{fig:bxpb}, the
relative improvement from the best PSWF to the best filter at all is very
slight (a few percents). Here, the most favorable feature of our approach is
that there is no tuning parameters (opening of the PSWF for instance, or the
order of the splines window functions if they are taken as alternatives) to be
found before the analysis.

\section{Conclusion}
\label{sec:conclusion}

A flexible way of analysing a field on the sphere in a space-frequency
manner has been presented. It is based on the needlet construction of
\cite{narcowich:petrushev:ward:2006}. The proposed analysis functions
form a frame in the space a square-integrable functions on the
sphere. Decompositions are essentially operating in the Spherical
Harmonics domain, leading to fast implementations. Various criteria
are used to design good spectral windows. This optimization can lead
to decisive improvement in high precision applications such as modern
cosmology (CMB spectral estimation, component separation, etc.), where
localized analysis is crucial.

\begin{ack}
  We wish to thank Jacques Delabrouille for fruitful discussions
  motivating this work for CMB analysis.  Numerical work was conducted
  using HEALPix~\citep{Gorski+2005}.
\end{ack}

\appendix

\section{Proofs}
\label{sec:proofs}

\paragraph*{Proof of Propositions~\ref{prop:tightframe} and~\ref{prop:frame}}

Propositions~\ref{prop:tightframe} is a particular case of
Proposition~\ref{prop:frame}. Indeed (\ref{reconstruction})-(\ref{eq:bsqrth})
imply (\ref{eq:boundssumblj}) with $C_1=C_2=1$. Together
with~(\ref{eq:defbtilde}) we get $\tilde\beta_k^{(j)} = \beta_k^{(j)}$ and
$\tilde \psi_k^{(j)} = \psi_k^{(j)}$. Prove now Proposition~\ref{prop:frame}.
Firstly, using successively (\ref{eq:kernel}) and the quadrature assumption
(remind that for any $\ell,\ell' \leq d$, $(\Pi_\ell X) (\Pi_\ell X) \in
\bigoplus_{l=0}^{2d}\Hset_{l}$)
\begin{align*}
  \sum_{j \in \mcJ, k \in K^{(j)}}|\beta_k^{(j)}|^2 &= \sum_{j\in\mcJ}\sum_{k \in K^{(j)}} \lambda_k^{(j)} \left|
  \sum_{\ell = 0}^{d^{(j)}}b_\ell^{(j)} \int X(\xi)L_\ell(\xi \cdot \xi_k) \dd \xi
\right|^2 \\ 
&= \sum_{j\in\mcJ}\sum_{k \in K^{(j)}} \lambda_k^{(j)} \left|
  \sum_{\ell = 0}^{d^{(j)}}b_\ell^{(j)} \Pi_\ell X (\xi_k) \right|^2 \\
&= \sum_{j\in\mcJ} \sum_{\ell, \ell' = 0}^{d^{(j)}}b_\ell^{(j)}b_{\ell'}^{(j)} 
\sum_{k\in K^{(j)}} \lambda_k^{(j)}
\Pi_\ell X (\xi_k) \Pi_{\ell'} X (\xi_k)\\
&= \sum_{j\in\mcJ} \sum_{\ell, \ell' = 0}^{d^{(j)}} b_\ell^{(j)}b_{\ell'}^{(j)} 
\int_\Sset \Pi_\ell X(\xi)\Pi_{\ell'} X(\xi) \dd \xi \\
&= \sum_{j\in\mcJ} \sum_{\ell, \ell' = 0}^{d^{(j)}} b_\ell^{(j)}b_{\ell'}^{(j)}  \delta_{\ell \ell'}
\int_\Sset | \Pi_\ell X(\xi)|^2 \dd \xi\\
&= \sum_{\ell \in \Nset}\sum_{j\in\mcJ} (b_\ell^{(j)})^2 \|\Pi_\ell X\|^2 \;
.
\end{align*}
Using (\ref{eq:boundssumblj}) and $\|X\|^2 = \sum_\ell \|\Pi_\ell X\|^2$, we get 
$
  C_1 \|X\|^2 \leq  \sum_{j,k}|\beta_k^{(j)}|^2  \leq   C_2 \|X\|^2
$. Prove now that $(\tilde \psi_k^{(j)})$ is the dual frame of
$(\psi_k^{(j)})$. Write
\begin{align*}
  \langle \tilde \psi_{k'}^{(j')} , \psi_{k}^{(j)} \rangle &=
  (\lambda_{k'}^{(j')}\lambda_{k}^{(j)})^{(1/2)}
  \sum_{\ell'=0}^{d^{(j')}}\sum_{\ell=0}^{d^{(j)}} \tilde b_{\ell'}^{(j')}
  b_{\ell}^{(j)} \int_\Sset L_{\ell'}(\xi \cdot \xi_{k'}^{(j')}) L_{\ell} (\xi
  \cdot \xi_{k}^{(j)}) \dd \xi \\
  &=   (\lambda_{k'}^{(j')}\lambda_{k}^{(j)})^{(1/2)} \sum_{\ell=0}^{d^{(j')}} \tilde b_{\ell}^{(j')} b_{\ell}^{(j)} L_{\ell}(\xi_{k'}^{(j')} \cdot \xi_{k}^{(j)})
\end{align*}
Then, for any $j \in \Nset$, $k\in K^{(j)}$
\begin{align*}
  \sum_{j',k'} \langle \tilde \psi_{k}^{(j)} , \psi_{k'}^{(j')} \rangle
  \psi_{k'}^{(j')} &=  (\lambda_{k}^{(j)})^{1/2} \sum_{j',k'} \lambda_{k'}^{(j')}  \sum_{\ell=0}^{d^{(j')}} \tilde b_{\ell}^{(j')} b_{\ell}^{(j)} L_{\ell}(\xi_{k'}^{(j')}
  \cdot \xi_{k}^{(j)}) \sum_{\ell'=0}^{d^{(j')}}b_{\ell'}^{(j')}
  L_{\ell'}(\xi_{k'}^{(j')}   \cdot \xi)\\
  &=  (\lambda_{k}^{(j)})^{1/2}  \sum_{j' \in \mcJ}\sum_{\ell=0}^{d^{(j')}} \tilde b_{\ell}^{(j')} b_{\ell}^{(j)}  \sum_{\ell'=0}^{d^{(j')}} b_{\ell'}^{(j')}
  \sum_{k' \in K^{(j)}}  \lambda_{k'}^{(j')} L_{\ell}(\xi_{k'}^{(j')}
  \cdot \xi_{k}^{(j)})   L_{\ell'}(\xi_{k'}^{(j')}   \cdot \xi) \\
&=  (\lambda_{k}^{(j)})^{1/2}  \sum_{j' \in \mcJ} \sum_{\ell=0}^{d^{(j')}} \tilde b_{\ell}^{(j')} b_{\ell}^{(j)}  \sum_{\ell'=0}^{d^{(j')}} b_{\ell'}^{(j')}
  \int_\Sset L_{\ell}(\xi'
 \cdot \xi_{k}^{(j)})   L_{\ell'}(\xi'   \cdot \xi) \dd \xi' \\
&=  (\lambda_{k}^{(j)})^{1/2}  \sum_{j' \in \mcJ} \sum_{\ell=0}^{d^{(j')}} \tilde b_{\ell}^{(j')} b_{\ell}^{(j)}  \sum_{\ell'=0}^{d^{(j')}} b_{\ell'}^{(j')}
  \delta_{\ell \ell'} L(\xi_k^{(j)} \cdot \xi) \\
&= (\lambda_{k}^{(j)})^{1/2}  \sum_{\ell=0}^{\infty}  b_{\ell}^{(j)}  \sum_{j'
  \in \mcJ}  \tilde b_{\ell}^{(j')}b_{\ell}^{(j')}  L(\xi_k^{(j)} \cdot \xi) \\
&= \psi_k^{(j)} \; .
 \end{align*}
The assertions~(\ref{eq:dualframes}) are a consequence of the dual frame
property \cite[see \eg{}][]{daubechies:1992}.

\paragraph*{Proof of Proposition \ref{prop:expneedcoeff}}
  From Definition~\ref{def:needlet} of $\psi_k$ and
  Eq~(\ref{eq:kernel})
  \begin{equation*}
    \beta_k  =  \langle X,\psi_k \rangle 
    = \sqrt{\lambda_k}\sum_\ell b_\ell\int_\Sset X(\xi)L_\ell(\xi,\xi_k)\dd\xi
    = \sqrt{\lambda_{k}}\Phi X(\xi_k) .   \qed  
  \end{equation*}

\paragraph*{Proof of Eq.~(\ref{eq:master})}
\begin{eqnarray*}
  (2\ell+1)\esp(\hat{C}_\ell)&=& \sum_{m=-\ell}^{\ell} \dint_{\Sset \times \Sset} \esp \left\{X(\xi)X(\xi')\right\}Y\lm(\xi)Y\lm(\xi')W(\xi)W(\xi')\dd\xi \dd\xi'\\
  &=& \dint_{\Sset \times \Sset} \bigl\{\sum_{\ell' \in  \Nset}C_{\ell'}L_{\ell'}(\xi \cdot \xi')\bigr\}L_\ell(\xi
  \cdot \xi')W(\xi)W(\xi')\dd\xi \dd\xi'\\
  &=&
  \sum_{\ell'\in\Nset}C_{\ell'}\sum_{0\leq\ell''\leq\ell+\ell'}\alpha_{\ell \ell' \ell''}\dint_{\Sset \times \Sset} 
  L_{\ell''}(\xi \cdot \xi')W(\xi)W(\xi')\dd\xi \dd\xi'\\
  &=& \sum_{\ell'\in\Nset}C_{\ell'}\sum_{0\leq\ell''\leq\ell+\ell'}\alpha_{\ell \ell' \ell''}(2\ell''+1)C_{\ell''}^W
\end{eqnarray*}

\paragraph*{Proof of Proposition~\ref{prop:expstatcriterion}}
 As $X\bar W = \sumstar_{(\ell_1,m_1)} a_{\ell_1 m_1}Y_{\ell_1 m_1}
\sumstar_{(\ell_2,m_2)} \bar w_{\ell_2 m_2}Y_{\ell_2 m_2}$, 
\begin{gather*}
  \bar a\lm := \langle X\bar W , Y\lm\rangle = \sumstar_{(\ell_i,m_i)_{i=1,2}} a_{\ell_1 m_1}\bar
  w_{\ell_2 m_2}\varwigner{\ell_1}{\ell_2}{\ell}{m_1}{m_2}{m} \ .
\end{gather*}
Together with  $\esp[a\lm a_{\ell' m'}^*] = C_\ell \delta_{\ell \ell'}\delta_{m
  m'}$ it yields
\begin{equation}
  \label{eq:covbaralm}
  \esp [\bar a\lm \bar a_{\ell' m'}^*]  = \sumstar_{(\ell_i,m_i)_{i=1,2,3}} C_{\ell_1} \bar
  w_{\ell_2 m_2} \bar w_{\ell_3 m_3}^*
 \varwigner{\ell_1}{\ell_2}{\ell}{m_1}{m_2}{m} \varwigner{\ell_1}{\ell_3}{\ell'}{m_1}{m_3}{m'}^*
\end{equation}
 Combining~(\ref{eq:expepsilon}) and~(\ref{eq:covbaralm}) we get
\begin{align*}
  R(\bb) &= ({\sum_{\ell\in\Nset} \sigma_\ell^2h_\ell^2})^{-1}\esp \int_\Sset
  \sumstar_{(\ell_4,m_4)} d_{\ell_4 m_4} Y_{\ell_4 m_4}(\xi)
  \left|\sumstar_{(l,m)} h_\ell \bar a\lm Y\lm(\xi)\right|^2 \dd \xi \\&=
 ({\sum_{\ell\in\Nset} \sigma_\ell^2h_\ell^2})^{-1}\sum_{\ell,\ell' \in
    \Nset}h_\ell h_{\ell'} \sum_{m m'} \esp [\bar a\lm \bar a_{\ell' m'}^*] 
 \sumstar_{(\ell_4,m_4)} d_{\ell_4 m_4} \int_\Sset Y _{\ell_4 m_4}(\xi)  Y\lm(\xi)
 Y_{\ell' m'}^*(\xi) \dd \xi
\\&=  ({\sum_{\ell\in\Nset} \sigma_\ell^2h_\ell^2})^{-1}\sum_{\ell,\ell' \in
    \Nset}h_\ell h_{\ell'} Q_{\ell \ell'} . 
\end{align*}
If $W$ is axisymmetric,
\begin{align}
  \esp [\bar a\lm \bar a_{\ell' m'}^*] &= \sumstar_{\ell_1,m_1} C_{\ell_1}
  \sum_{\ell_2,\ell_3} \bar w_{\ell_2 0} \bar w_{\ell_3 0}
  \varwigner{\ell_1}{\ell_2}{\ell}{m_1}{0}{m}
  \varwigner{\ell_1}{\ell_3}{\ell'}{m_1}{0}{m'} \nonumber \\
  &= \delta_{m,m'} \sum_{\ell_1,\ell_2,\ell_3} C_{\ell_1} \bar w_{\ell_2 0}
  \bar w_{\ell_3 0} \varwigner{\ell_1}{\ell_2}{\ell}{m}{0}{m}
  \varwigner{\ell_1}{\ell_3}{\ell'}{m}{0}{m} =: A_{\ell \ell' m}
 \label{eq:covbaralmaxisym}
\end{align}
where we used the fact that $\varwigner{\ell}{\ell'}{\ell''}{m}{0}{m''} = 0$ if
$m \neq m''$ and that $w_{l_3 0}$ and $ \varwigner{\ell_1}{\ell_2}{\ell}{m}{0}{m}$
are real. 
If $D$ is axisymmetric and with $ D_{\ell \ell' m} := \sum_{\ell_4}d_{\ell_4,0}
\varwigner{\ell}{\ell_4}{\ell'}{m}{0}{m}$
\begin{align*}
  Q_{\ell \ell'} &=
  \sum_{m}\sum_{(\ell_i,m_i)_{i=1,2,3}}C_{\ell_1}
   \bar w_{\ell_2 m_2} \bar w_{\ell_3 m_3}^*
  \varwigner{\ell_1}{\ell_2}{\ell}{m_1}{m_2}{m}
  \varwigner{\ell_1}{\ell_3}{\ell'}{m_1}{m_2}{m}^*  \sum_{\ell_4
    \in\Nset}d_{\ell_4 0} \varwigner{\ell}{\ell_4}{\ell'}{m}{0}{m} . \;  \qed
\end{align*}

\section{Legendre polynomials, spherical harmonics and related useful formulae}
\label{sec:leg-pol}

Usually, $P_\ell(z)$ denotes the Legendre polynomial of order $\ell$,
normalized by $P_\ell(1)=1$.  For our purposes, it is more convenient to use a
different normalization
\begin{displaymath}
  L_\ell (z) = \frac{2\ell+1}{4\pi} P_\ell (z)
\end{displaymath}
because we get coefficient-free properties like
\begin{displaymath}
  L_\ell ( \xi' \cdot \xi ) =
  \sum_{m=-\ell}^{\ell} Y_{\ell m}^*(\xi) Y_{\ell m} (\xi ')
\end{displaymath}
and
\begin{equation}
\label{eq:intLellLellp}
\int_{\Sset} L_\ell    (\eta \cdot \xi) L_{\ell'} (\eta' \cdot \xi )
\dd\xi = \delta_{\ell\ell'} L_\ell(\eta \cdot  \eta').
\end{equation}
In other words, $L_\ell$ is the polynomial kernel of the harmonic
projection on $\mathbb{H}_\ell$. 
We have $\int_{-1}^{+1} P_\ell(z)^2 \dd z = \frac2{2\ell+1}$ and  $\int_{-1}^{+1} L_\ell(z)^2 \dd
z = \frac{2\ell+1}{8\pi^2}$.

The spherical harmonics are explicitly given in a factorized form in
terms of the associated Legendre polynomials and the complex
exponentials as
\begin{displaymath}
  Y\lm(\theta,\varphi)=\sqrt{\frac{(2\ell+1)}{4\pi}\frac{(\ell-m)!}{(\ell+m)!}} P\lm(\cos\theta)e^{im\varphi}
\end{displaymath}
where $P\lm(x)=(-1)^m(1-x^2)^{m/2}\frac{\mathsf{d}^m}{\mathsf{d}x^m}P_\ell(x)$.

The following equations relate the integral of the product of three
complex spherical harmonics over the total solid angle or three
Legendre polynomials with the Wigner-3$j$ coefficients (for a
definition in terms of Clebsh-Gordan coefficients, see
\cite{varshalovich:etal:1988}, pp235--).
\begin{align}
  \label{eq:wigner3j1}
  \varwigner{\ell_1}{\ell_2}{\ell_3}{m_1}{m_2}{m_3} &= 
  \int_\Sset Y_{\ell_1 m_1}(\xi)Y_{\ell_2 m_2}(\xi)Y^*_{\ell_3 m_3}(\xi) \dd
  \xi \nonumber \\ &= (-1)^{m_3}
\int_\Sset Y_{\ell_1 m_1}(\xi)Y_{\ell_2 m_2}(\xi)Y_{\ell_3 -m_3}(\xi) \dd \xi
\nonumber \\ &= 
(-1)^{m_3}  \sqrt{\frac{(2\ell_1+1)(2\ell_2+1)(2\ell_3+1)}{4\pi}}
  \wigner{\ell_1}{\ell_2}{\ell_3}{0}{0}{0}\wigner{\ell_1}{\ell_2}{\ell_3}{m_1}{m_2}{-m_3}
\end{align}
\begin{align}
\label{eq:wigner3j2}
  \frac 12 \int L_\ell(z)
  L_{\ell'}(z')L_{\ell''}(z'') \dd z \dd z'\dd z'' &= \frac{(2\ell+1)(2\ell'+1)(2\ell''+1)}{(4\pi)^3}
  \wigner{\ell}{\ell'}{\ell''}{0}{0}{0}^2 \\
  &= (4\pi)^{-2} \varwigner{\ell_1}{\ell_2}{\ell_3}{0}{0}{0}
\end{align}
From (\ref{eq:wigner3j2}) and $\int L_\ell
  L_{\ell'}=\delta_{\ell,\ell'}\frac{2\ell+1}{8\pi^2}$ we get:
$$L_\ell L_{\ell'}=\sum\limits_{0\leq\ell''\leq\ell+\ell'}\alpha_{\ell \ell' \ell''}L_{\ell''}$$ with
\begin{equation}
\label{eq:defalpha}\alpha_{\ell \ell' \ell''}=\frac{(2\ell+1)(2\ell'+1)}{4\pi}\wigner{\ell}{\ell'}{\ell''}{0}{0}{0}^2
\ .
\end{equation}

\bibliography{waveletdesign}
\bibliographystyle{elsart-harv}

\end{document}